\documentclass[a4paper,11pt]{article}
\usepackage{amssymb,amsfonts,latexsym,amsmath,amsthm,amscd,graphpap,color}

\usepackage{graphicx}

\usepackage{fancybox,framed,xcolor,dsfont}
\newtheorem{theorem}{Theorem}[section]
\newtheorem{prop}[theorem]{Proposition}
\newtheorem{lemma}[theorem]{Lemma}
\newtheorem{corollary}[theorem]{Corollary}
\newenvironment{demo}{ \noindent \emph{\textbf{Proof:}}}{\hfill$\square$\\
\vspace{0.4cm}}

\newcommand{\RR}{\mathbb{R}}

\newcommand{\Lc}{\mathcal{L}}

\newcommand{\NN}{\mathbb{N}}

\newcommand{\Cc}{\mathcal{C}}

\newcommand{\Mc}{\mathcal{M}}

\newcommand{\Oc}{\mathcal{O}}

\newcommand{\Ac}{\mathcal{A}}

\newcommand{\ZZ}{\mathbb{Z}}

\newcommand{\Un}{1\hspace{-1.5mm}1}
\renewcommand{\Re}{\mathrm{Re}}

\newcommand{\Op}{{\mathrm{Op}}}
\newcommand{\supp}{\mathrm{supp}}
\newcommand{\no}{n$^{\text{o}}$}
\newcommand{\grad}{\nabla}
\renewcommand{\div}{\mathrm{div}}

\renewcommand{\t}{^{\intercal}}
\newcommand{\ul}{{\text{ul}}}
\newcommand{\loc}{{\text{loc}}}

\newcommand{\pc}{ \usefont{T1}{cmtl}{m}{n} \selectfont}
\newcounter{compte}
\newenvironment{enum2}{\begin{list}{\roman{compte})} {\usecounter{compte}
\topsep=1mm \itemsep=0.2mm \leftmargin=5mm  } }
{\end{list}}
\numberwithin{equation}{section}

\newdimen\texpscorrection
\newdimen\figcenter
\def\figurewithtex #1 #2 #3 #4 #5\cr{\null
  {\goodbreak\figcenter=\hsize\relax
  \advance\figcenter by -#4truecm
  \divide\figcenter by 2
  \begin{figure}[hbt]
  \vskip #3truecm\noindent\hskip\figcenter
  \includegraphics{#1}{\hskip\texpscorrection\input #2 }
  \vskip 0.8truecm{\baselineskip=0.8\baselineskip
  \noindent \vbox{\noindent {\footnotesize #5}}\par}
  \end{figure}}}
\def\point#1 #2 #3 {\rlap{\kern #1 truecm
\raise #2 truecm \hbox{#3}}}

\colorlet{shadecolor}{white!93!black}

\begin{document}

\title{\bf Exponential decay for the damped wave equation in unbounded domains}

\author{Nicolas \textsc{Burq}\footnote{D\'epartement de Math\'ematiques d'Orsay
- UMR 8628 CNRS/Universit\'e Paris-Sud - Bat. 425 - F-91405 Orsay Cedex, France
email: {\pc nicolas.burq@math.u-psud.fr}}~\&~Romain
\textsc{Joly}\footnote{Institut Fourier - UMR5582
CNRS/Universit\'e de Grenoble - 100, rue des Maths - BP 74 - F-38402
St-Martin-d'H\`eres, France, email: {\pc romain.joly@ujf-grenoble.fr}}}

\date{}

\maketitle

\begin{abstract}
We study  the  decay of the semigroup generated by the damped wave
equation in an unbounded domain. We first prove under the natural {\em geometric
control condition} the exponential decay of the semigroup. Then we prove under
a weaker condition the logarithmic decay of the solutions (assuming that the
initial data are smoother). As corollaries, we obtain several extensions
of previous results of stabilisation and control.

\vspace{3mm}

On \'etudie la d\'ecroissance du semi-groupe des ondes amorties dans un domaine
non born\'e. Notre premier r\'esultat est que, sous une hypoth\`ese naturelle de
contr\^ole g\'eom\'etrique, le semi-groupe d\'ecro\^\i{}t exponentiellement
vite. On d\'emontre ensuite sous une hypoth\`ese plus faible la d\'ecroissance
logarithmique des solutions associ\'ees \`a des donn\'ees initiales plus
r\'eguli\`eres. On obtient en corollaire plusieurs g\'en\'eralisations de
r\'esultats de stabilisation et de contr\^ole.

\vspace{3mm}

\end{abstract}

%%%%%%%%%%%%%%%%%%%%%%%%%%%%%%%%%%%%%%%%%%%%
%%%%%%%%%%%%%%%%%%%%%%%%%%%%%%%%%%%%%%%%%%%%

\section{Introduction}

In this article we consider the damped wave equation. In the simplest case of
constant coefficients Laplace operator, our main result takes the following
form: 
\begin{theorem}
Let $\gamma\in L^\infty(\RR^d)$ be a non-negative damping. Assume that
$\gamma$ is a uniformly continuous function and that there exist $L,c>0$ such
that for any $(x_0, \xi_0) \in \mathbb{R}^d\times \mathbb{S}^{d-1} $,
$$ \int_{s=0}^L \gamma(x_0 + s \xi_0)dx~\geq~c~>~0~.$$
Then, there exist $M$ and $\lambda>0$ such that any solution of
$$\partial^2_{tt} u +\gamma(x) \partial_t u= (\Delta-Id) u
\hphantom{espace} (t,x)\in\RR_+\times\RR^d$$
satisfies 
$$\|u(t)\|_{H^1(\RR^d)}+\|\partial_tu(t)\|_{L^2(\RR^d)}~\leq~Me^{
-\lambda t}\left( \|u(0)\|_{H^1(\RR^d)}+\|\partial_tu(0)\|_{L^2(\RR^d)}
\right)~.$$
\end{theorem}

\subsection{The damped wave equation:}
More precisely, our results concern a more general linear damped wave equation
in $\RR^d$, with $d\geq 1$:
\begin{equation}\label{eq}
\left\{\begin{array}{ll}
\partial^2_{tt} u(x,t) +\gamma(x) \partial_t u(x,t)=
\div(K(x)\grad u(x,t))-u(x,t) ~~~~&(t,x)\in\RR_+\times\RR^d~,\\
(u,\partial_t u)(\cdot,0)=(u_0,u_1)\in H^1(\RR^d)\times L^2(\RR^d)&
\end{array}\right.
\end{equation}
where $K\in \Cc^\infty(\RR^d,\Mc_{d}(\RR))$ is a smooth family of real symmetric
matrices, which are uniformly positive in the sense that there exist
two positive constants $K_{\rm inf}$ and  $K_{\rm sup}$ such that
\begin{equation}\label{condi-K}
\forall \xi\in\RR^d~,~~K_{\rm sup}|\xi|^2\geq \xi\t.K(x).\xi \geq
K_{\rm inf}|\xi|^2~. 
\end{equation}
The damping coefficient $\gamma\in L^{\infty}(\RR^d)$ is assumed to be a
bounded and non-negative function. We set $X=H^1(\RR^d)\times L^2(\RR^d)$ and 
\begin{equation}\label{def-A}
A=\left(\begin{array}{cc} 0& Id\\(\div(K(x)\grad)-Id)&
-\gamma(x)\end{array}\right)\hphantom{espace}D(A)=H^2(\RR^d)\times
H^1(\RR^d)~.
\end{equation}
We equipped $H^1(\RR^d)$ with the scalar product
\begin{equation}\label{prod-scalaire-H1}
\langle u|v\rangle_{H^1}=\int_{\RR^d} (\grad u(x))\t.K(x).(\overline {\grad
v(x)})\,+\,u(x)\overline{v(x)}\,dx~. 
\end{equation}
Obviously, this scalar product is equivalent to the classical one and direct
computations show that it satisfies 
$$\langle (\div(K(x)\grad)-Id)u|v\rangle_{L^2}=-\langle u|v\rangle_{H^1}
~~~\text{ and }~~~\Re(\langle AU|U\rangle_X)=-\int \gamma(x)|v(x)|^2\,dx~
$$
for any $U=(u,v)\in D(A)$. Then, one easily checks that $A$ is a dissipative
operator and therefore generates a semigroup $e^{At}$ on $X$.
\subsection{ Exponential decay and Hamiltonian flow:}
The main purpose of this paper is to investigate the exponential decay of the
semigroup associated to \eqref{eq}: we ask whether there exist $M$ and $\lambda>0$
such that
\begin{equation}\label{exp-dec-1}
\forall t\geq0~,~~|||e^{At}|||_{\Lc(X)}\leq Me^{-\lambda t}~.
\end{equation}
For the damped wave equation in a bounded domain and a continuous damping coefficient, it is well
known that the exponential decay is equivalent to the fact that all the
trajectories of the Hamiltonian flow intersect the support of the damping (see
\cite{RT}, \cite{BLRinterne}, \cite{BLR}) and~\cite{BuGe}. More precisely, to the Laplacian
operator with variable coefficients $\div(K(x)\grad)$, we associate the symbol
$g(x,\xi)=\xi\t.K(x).\xi$ and the Hamiltonian flow
$\varphi_t(x_0,\xi_0)=(x(t),\xi(t))$ defined on $\RR^{2d}$ by
\begin{equation}\label{eq-hamiltonien}
\varphi_0(x_0,\xi_0)=(x_0,\xi_0)~~~~\text{ and }~~~~\partial_t
\varphi_t(x,\xi)=(\partial_\xi g(x(t),\xi(t)),-\partial_x
g(x(t),\xi(t))~.
\end{equation}
We introduce the mean value of the damping along a ray a length $T$:
\begin{equation}\label{def-moy-gamma}
\langle \gamma \rangle_T(x,\xi)=\frac1T\int_0^T \gamma(\varphi_t(x,\xi))dt
\end{equation}
where we use the obvious notation $\gamma(x,\xi):=\gamma(x)$. We also introduce
the set $\Sigma$ of rays of speed one, that is 
\begin{equation}\label{def-sigma}
\Sigma=\{(x,\xi)\in\RR^{2d}~,~\xi\t K(x) \xi=1\}~.
\end{equation}

\vspace{3mm}

{\noindent \bf Some previous works:}\\
If $\Omega$ is a bounded manifold, the uniform positivity of $\langle \gamma
\rangle_T(x,\xi)$ in $\Sigma$ for some $T>0$ implies that the exponential decay
\eqref{exp-dec-1} holds, as shown in the celebrated articles \cite{RT},
\cite{BLRinterne} and \cite{BLR} of Bardos, Lebeau, Rauch and Taylor. The
assumption that there exists $T>0$ such that $\langle \gamma
\rangle_T(x,\xi)>0$ in $\Sigma$ is called the {\it geometric control condition}.
The article \cite{Lebeau} underlines in addition the importance of the value of
$\min_{(x,\xi)\in\Sigma}\,\langle \gamma \rangle_T(x,\xi)$ in order to control
the rate of decay of the high frequencies.

In the case of an unbounded manifold, two situations have been investigated.
First, some authors have considered the free wave equation \eqref{eq} in an
exterior domain (with $\gamma\equiv 0$ or $\gamma>0$ only on a compact subset
the exterior domain). They have shown that
the local energy decays to zero in the sense that, under suitable assumptions,
the energy of any solution escapes away from any compact set, see \cite{LMP},
\cite{MRS} and \cite{AK} and the references therein. Secondly, several works
have studied the damped wave equation in an unbounded manifold and with a
non-linearity, but assuming that the damping satisfies $\gamma(x)\geq \alpha>0$
outside a compact set, see \cite{Zuazua2}, \cite{Feireisl}, \cite{DLZ} and
\cite{Jol-Lau}. 

Considering these previous works, it appears that one natural case has not been
studied: the exponential decay of the semigroup $e^{At}$ generated by the damped
wave equation on a whole unbounded manifold, with the geometric control
condition only, that is without assuming that $\gamma\geq \alpha>0$ outside a
compact set. To our knowledge, this case is surprisingly missing in the
literature. The main purpose of this article is to settle this natural problem.

\vspace{3mm}

{\noindent \bf Main results:}\\
We denote by $\Cc^k_b(\RR^d)$ the set of functions in $\Cc^k(\RR^d)$ which are
bounded, as well as their $k$ first derivatives. If $k=\infty$, the bound is
not assumed to be uniform with respect to the derivatives. 
We recall that $\langle \gamma \rangle_T$ and $\Sigma$ have been defined
in \eqref{def-moy-gamma} and \eqref{def-sigma}. Our main
result is as follows.
\begin{theorem}\label{th-main}
\begin{shaded}
We assume that the metric $K$ belongs to
$\Cc^\infty_b(\RR^d,\Mc_d(\RR))$ and that the bounded non-negative damping
$\gamma$ is uniformly continuous and satisfies 
$$(GCC)~\hphantom{espace}~\text{ there exist }T,\alpha>0\text{ such that
}\langle
\gamma \rangle_T(x,\xi)\geq\alpha>0~,~\text{ for all }
(x,\xi)\in\Sigma~.$$

Then, the semigroup generated by the damped wave equation \eqref{eq} is
exponentially decreasing that is that there exist $M$ and $\lambda>0$ such
that
\begin{equation}\label{exp-dec}
\forall t\geq 0~,~~|||e^{At}|||_{\Lc(X)}\leq Me^{-\lambda t}~.
\end{equation}
\end{shaded}
\end{theorem}

Assume now that the geometric control condition (GCC) is violated but the
damping is still efficient on a network of balls. The Lasalle invariance
principle ensures that, for any initial data, the energy of the solution goes to
$0$ when $t \rightarrow + \infty$. Since the geometric control condition does
not hold, it is classical that the convergence to $0$ can be arbitrarily slow:
$$ \forall T>0, \sup_{(u_0, u_1)\in H^1\times L^2 }\frac{ \| (u(T), \partial_t
u(T))\|_{H^1 \times L^2} } {\|(u_0, u_1)\|_{H^1\times L^2}} =1. $$
Our second result extends to the non compact setting a similar result by Lebeau
proved on compact manifolds~\cite{Lebeau} (see also~\cite{LR}) and 
gives an upper bound for the rate of decay when the initial data are smoother
(see~\cite[Definition 1.1 and Section 3.1]{LRM} for a similar geometric setting
developped independently).
\begin{theorem} \label{Th.2}
\begin{shaded}
We assume that the metric $K$ belongs to $\Cc^\infty_b(\RR^d,\Mc_d(\RR))$
and that $\gamma\in L^\infty(\RR^d)$ satisfies 
\begin{align*}
 (NCC)~\hphantom{esp}~&\text{ there exist }L,r,a>0 \text{ and a sequence
}(x_n)\subset\RR^d~\text{ such that }\\
&\gamma(x)\geq a>0\text{ on }\cup_n
B(x_n,r)~~\text{ and }\forall x\in\RR^d~,~ d(x,\cup_n \{x_n\})\leq L~.
\end{align*}
Then, for any $k>0$, there exists $C_k>0$ such that for any $(u_0, u_1)\in
H^{k+1}( \mathbb{R}^d) \times H^k ( \mathbb{R}^d) $,
$$ \| (u(T), \partial_t u(T))\|_{H^1 \times L^2}  \leq \frac{ C_k } { \log (2+t)
^k} \|(u_0, u_1)\|_{H^{k+1} \times H^k}.$$
\end{shaded}
\end{theorem}

\vspace{3mm}

{\noindent \bf Some extensions and applications:}
\begin{enum2}
\item A contradiction argument shows very easily that, as soon as the
exponential decay holds for a damping coefficient $0\leq \gamma$, it holds (with
different constants, possibly worse) for any damping coefficient
$\widetilde{\gamma} \in L^\infty( \mathbb{R}^d)$ satisfying
$\tilde\gamma\geq\gamma$ (see the arguments of the second step of Section
\ref{sect-proof}). Consequently, Theorem \ref{th-main} also holds for any
damping $\gamma\in L^\infty(\RR^d)$ for which there exists $\underline{\gamma}$
with $0\leq\underline \gamma\leq \gamma$ satisfying (GCC) and being uniformly
continuous. Notice that the existence of $\underline{\gamma}$ uniformly
continuous satisfying (NCC) and $0\leq\underline \gamma\leq \gamma$ is
automatic in the case of Theorem \ref{Th.2}. That is why, the uniform
continuity can be omitted in its statement.

\item Theorem \ref{th-main} concerns solutions of \eqref{eq} with finite
energy. It is possible to consider solutions of \eqref{eq} with infinite energy
in the framework of uniformly local Sobolev spaces. The stabilisation in this
case is a straightforward corollary of Theorem \ref{th-main}, see Section
\ref{sect-applis}.

\item The ideas of the proof of Theorems \ref{th-main} and \ref{Th.2} may
apply to other geometric situations. For example, if we consider an unbounded
manifold
without boundary as a cylinder instead of $\RR^d$, then Theorems
\ref{th-main} and \ref{Th.2} will also hold with the obvious modifications of
their statements.

\item The smoothness assumptions on the coefficients $K(x)$ could be relaxed
(probably up to $\Cc^2$, see~\cite{Bu}). To keep the paper short, we chose not
to develop this issue here.

\item The exponential decay of the linear semigroup has important
applications in the control theory and the study of dynamics for the wave
equations. Some new results are obtained as corollaries of Theorem \ref{th-main}
as explained in Section \ref{sect-applis}. 
\end{enum2}

\vspace{3mm}

{\noindent \bf Remarks:}
\begin{enum2}
\item The simplest applications of Theorem \ref{th-main} are the
periodic frameworks satisfying the geometric control condition, see for example
Figure 1.a). To our knowledge, the exponential decay of the semigroup was not
known in this simple case (notice that one cannot directly use the framework of
the torus since the initial data $(u_0,u_1)$ are not periodic).

\item The proof of Theorem \ref{th-main} follows the lines of the
proofs of the results on compact manifolds (see \cite{BLRinterne},
\cite{BLR}, \cite{Zworski}\ldots). It also uses classical properties of
pseudo-differential calculus (see e.g. \cite{Alinhac}, \cite{Martinez} or
\cite{Lerousseau}). The main point in the analysis is to be careful when using
the classical arguments to deal with the infinity (in space). In particular this
forbids the use of tools as the defect measure, which only yields informations
on a compact subset of the domain. As usual, the proof of the stabilisation
stated in Theorem \ref{th-main} splits into two parts. The first is the control
of the high frequencies, where we fully use the geometric control condition
(GCC). This part is contained in Section \ref{sect-high}, where we have to
return to the semiclassical analysis behind the classical defect measure
arguments. The second part is the control of the low frequencies by using a
Carleman estimate as shown in Section \ref{sect-low}. In this section, we do not
use (GCC) but the weaker hypothesis (NCC), which is a uniform control of the
damping on a network of balls stated in Theorem \ref{Th.2}.

\item In a first version of this article by the second author alone, it
was shown that Theorem \ref{th-main} can be obtained in dimension
one by multipliers techniques following the ideas of \cite{Lions}
and \cite{Zuazua}. In some simple geometrical situation in higher dimension,
the multipliers techniques should also apply. The interest of this kind of
proofs is to provide explicit constants $M$ and $\lambda$, but the geometrical
assumptions cannot be as general as the ones of the main result of this paper,
except in dimension one.

\item Of course, our theorem also hold when the
operator $\div(K(x)\grad)-Id$ is replaced by the operator $\div(K(x)\grad)-\varepsilon Id$
with $\varepsilon>0$. However, when $\varepsilon=0$, that is when the
right-hand-side is not a negative operator but only a non-positive one, it is
known that one cannot hope an exponential decay of the solutions. Indeed, it
has been established since a long time (see \cite{Matsumura}) that the solutions
of $\partial_{tt}^2u+\partial_t u = \Delta u$ in $\RR^d$ asymptotically behave
as the ones of $\partial_t v = \Delta v$ (see for example \cite{OPZ}, \cite{RTY}
and the references therein). It is shown in \cite{CH} that if $u$ is
solution of $\partial_{ tt}^2u+\partial_t u = \Delta u$ in $\RR^d$, with initial
data $(u_0,u_1)\in H^1(\RR^d)\times L^2(\RR^d)$ and $v$ is solution of
$\partial_t v = \Delta v$ in $\RR^d$, with initial data $u_0+u_1$, then
$\|(u-v)(t)\|_{H^1}\leq C/t$. In particular, $u$ is generally decaying not
faster than $C/t^{d/4}$ for $d\leq 3$.

\item The uniform continuity assumption on $\gamma$ in Theorem
\ref{th-main} ensures that it can be regularised into a smooth
damping coefficient $\underline\gamma$ satisfying $\underline\gamma\leq\gamma$
and belonging to $\Cc^\infty_b(\RR^d,\RR)$. In particular, the fact that the
derivative of $\underline\gamma$ can be taken uniformly bounded will be
important in our proof order to apply the pseudo-differential calculus (notice
that these uniform
bounds would also be required if we used the multipliers techniques, at least
for the first derivative).
Of course, in the usual compact case, this assumption is
automatically satisfied. In Figures 1.b) and 1.c), we show examples where all
the Hypotheses of Theorem \ref{th-main} apply, if one neglects the
regularity hypothesis. In these cases, it would be natural to expect the
exponential decay of
the semigroup, but this is still an open problem. Notice that the simple
requirement that $\gamma$ belongs to $L^\infty$ is not sufficient to define
properly the mean value $\langle \gamma\rangle_T(x,\xi)$ everywhere. This could
be a hint that the regularisation assumption is not just a technical one.
\end{enum2}

\vspace{3mm}

{\noindent \bf Acknowledgements:} The authors thank Jean-Fran\c{c}ois Bony, 
Yves Colin de Verdi\`ere, Julien Royer and Patrick G\'erard for fruitful
discussions. The first author is partially supported by the Agence Nationale de
la Recherche through ANR-13-BS01-0010-03 (ANA\'E) and ANR-201-BS01-019-01
(NOSEVOL).

\begin{figure}[ht]
\ovalbox{
\begin{minipage}{16cm}
\begin{minipage}{7.5cm}
\begin{center}
\includegraphics[width= 7.5cm]{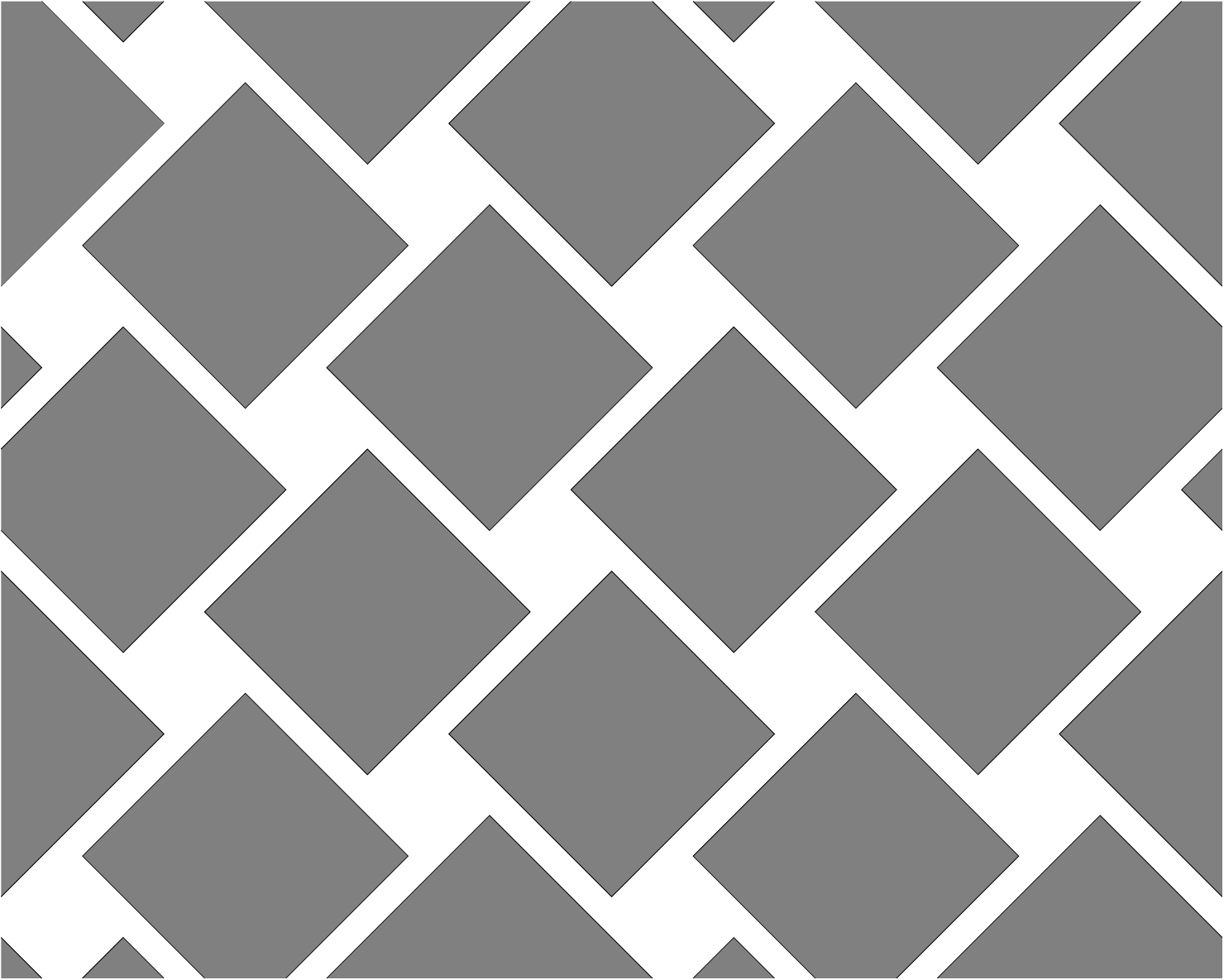}\\[2mm]
Figure 1.a): {A periodic two-dimensional example for which Theorem
\ref{th-main} holds: the semigroup generated by the corresponding damped wave
equation is decaying exponentially fast.\hfill}
\end{center}
\end{minipage}
~\vrule~
\begin{minipage}{8cm}
\begin{center}
\includegraphics[width= 8cm]{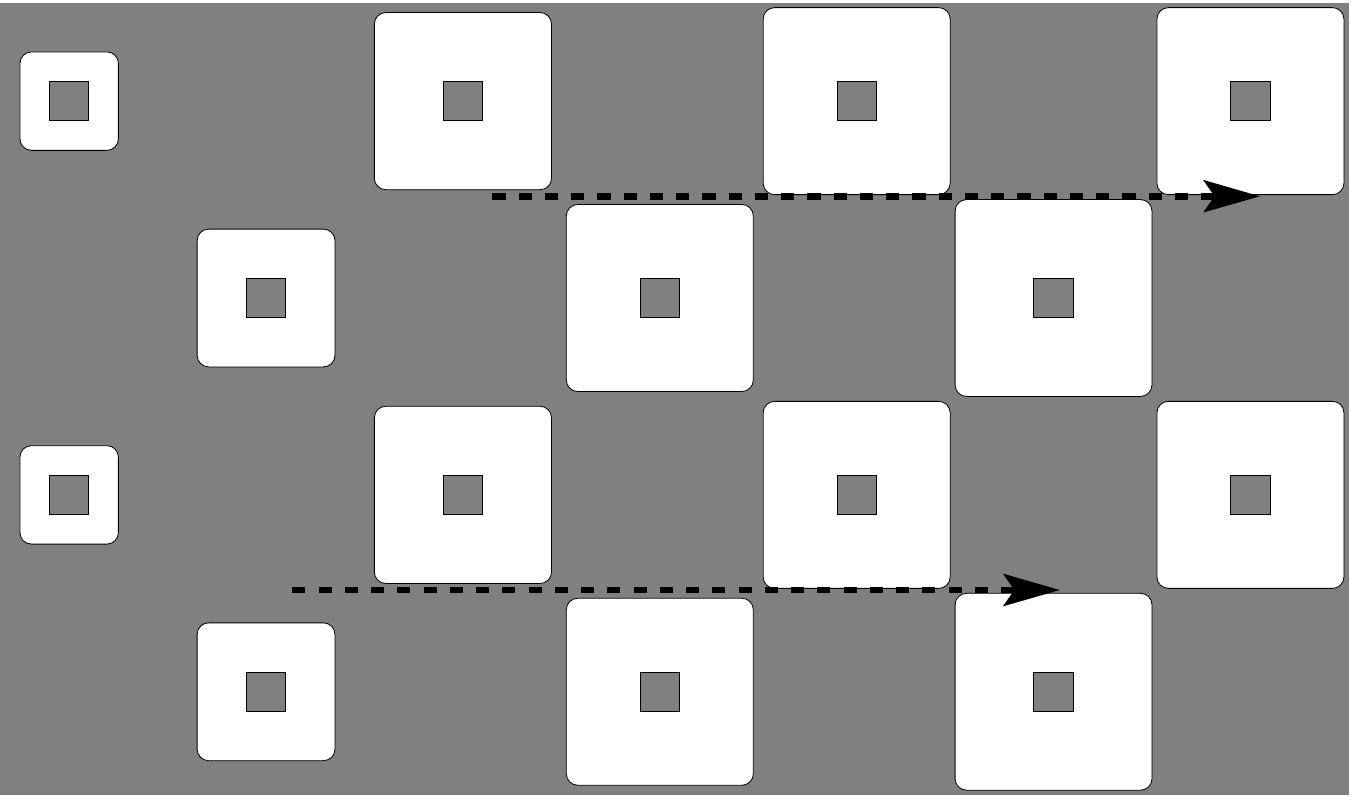}\\[2mm]
Figure 1.b): {A two-dimensional quasi-periodic exam[le where only the
regularisation condition in Hypothesis i) is not satisfied. Theorem
\ref{th-main} fails to apply because for any uniformly
continuous damping $\tilde\gamma$ satisfying $\tilde\gamma\leq \gamma$, the
infimum of the mean value $\langle \tilde\gamma\rangle$ is equal to $0$.\hfill}
\end{center}
\end{minipage}

\vspace{2mm}

\hrule

\vspace{3mm}
\begin{center}
\includegraphics[width= 12cm]{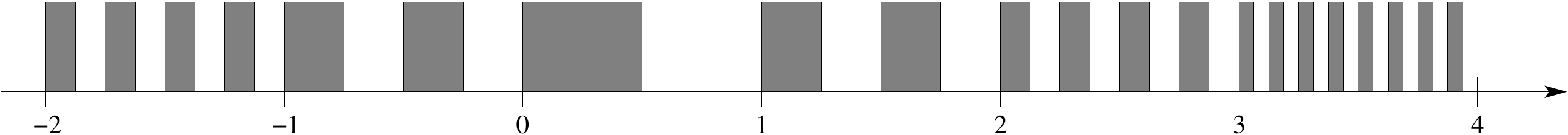}\\
Figure 1.c): {A one-dimensional example where the mean value of the damping
$\langle \gamma\rangle_1(x,\xi)$ is uniformly positive in $\Sigma$, but where
Theorem \ref{th-main} does not apply since there is no uniformly
continuous regularisation
$\tilde\gamma$ with mean value $\langle \tilde\gamma\rangle_T$
uniformly positive for some $T$ and with a derivative uniformly bounded in
$\RR$.\hfill}
\end{center}
\vspace{2mm}

\hrule

\vspace{0mm}
\begin{center}
\includegraphics[width=13cm]{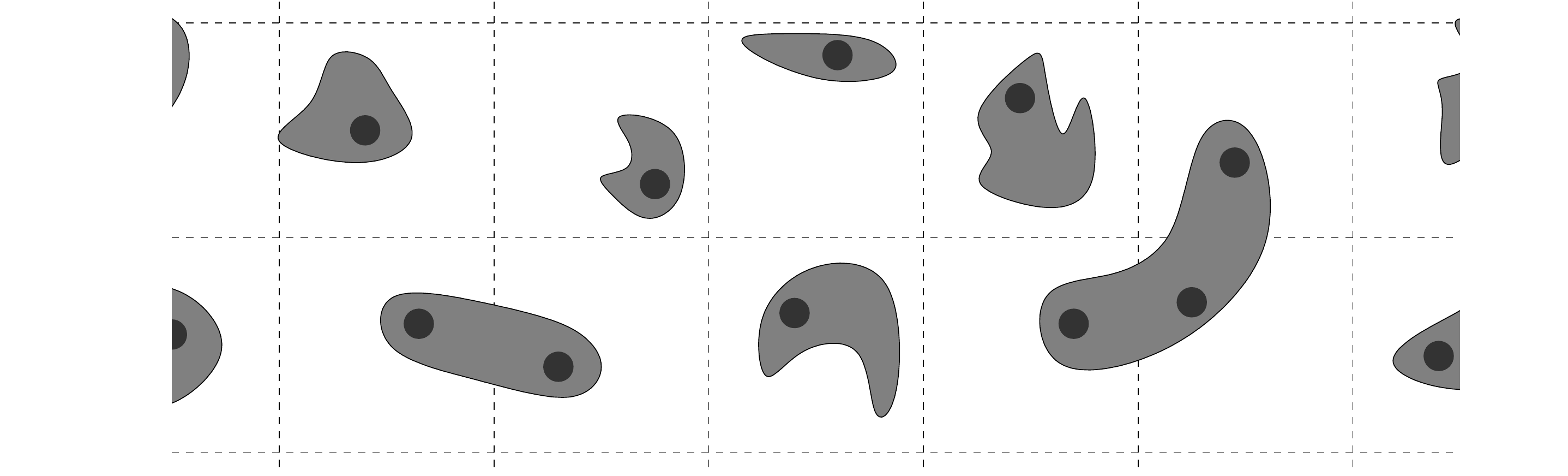}\\
 Figure 1.d): {An example where (NCC) holds but (GCC) does not. A network
of balls where the damping is effective is in dark grey. In this case, the
exponential decay of Theorem \ref{th-main} fails but the logarithmic decay of
Theorem \ref{Th.2} holds.\hfill}

\vspace{2mm}

\hrule

\vspace{0mm}

\caption{\it Discussion on some examples of damping. In the two-dimensional
situations, the damping is equal to $1$ on the grey regions and equal to $0$ in
the other regions. In the one-dimensional situation, the figure represents the
graph of the damping. In both cases, the metric is assumed to
be flat, i.e. $K(x)=Id$.}
\end{center}
\label{fig1}
\end{minipage}}
\end{figure}

%%%%%%%%%%%%%%%%%%%%%%%%%%%%%%%%%%%%%%%%%%%%
%%%%%%%%%%%%%%%%%%%%%%%%%%%%%%%%%%%%%%%%%%%%

\section{Proof of Theorem \ref{th-main}}\label{sect-proof}

In this section, we outline the proof of Theorem \ref{th-main}. The real
technical parts of its proof will be detailed in Sections \ref{sect-high} and
\ref{sect-low}.

There exist several ways to obtain the exponential decay \eqref{exp-dec} of the
semigroup $e^{At}$. The most classical one is to argue by contradiction to
establish the observation inequality $E(v(0))\leq C \int_0^T \gamma |\partial_t
v|^2dt$ for some $T>0$ and any solution $v$ of the free wave equation (see for
example \cite{Haraux} for the relation between this observation estimate and the
exponential decay of the damped semigroup). A less usual method consists in
uniformly estimating the resolvent $(A-\lambda Id)^{-1}$ on the imaginary axis
(see for example chapter 5 of \cite{Zworski}). We use here this last method as
a direct corollary of the results of \cite{Ge}, \cite{Pruss}
and \cite{Huang}. 

\vspace{3mm}

\noindent {\it $\bullet$ First step: a characterisation in terms of resolvent estimates.}\\
To study the exponential decay, we use here the characterisation given by
Theorem 3 of \cite{Huang}. 
\begin{theorem}[Gearhart-Pr\"uss-Huang]\label{th-Huang}
Let $e^{At}$ be a $\Cc^0-$semigroup in a Hilbert space $X$ and assume that
there exists a positive constant $M>0$ such that $|||e^{At}|||\leq M$ for
all $t\geq 0$. Then $e^{At}$ is exponentially stable if and only if
$i\RR\subset \rho(A)$ and 
\begin{equation}\label{eq-Huang}
\sup_{\mu\in\RR} |||(A-i\mu Id)^{-1}|||_{\Lc(X)}~<~+\infty~.
\end{equation}
\end{theorem}

Since the linear operator $A$ associated to the damped wave equation is
dissipative, we have $|||e^{At}|||\leq 1$ for all $t\geq 0$. To prove Theorem
\ref{th-main}, it remains to show that \eqref{eq-Huang} holds. We argue by
contradiction and assume that there exist two sequences
$(U_n)=(u_n,v_n)\subset D(A)=H^2(\RR^d)\times H^1(\RR^d)$ and
$(\mu_n)\subset\RR$ such that 
\begin{equation}\label{eq-absurde}
\|U_n\|_X^2=\|u_n\|_{H^1}^2+\|v_n\|_{L^2}^2=1~~~\text{ and
}~~~(A-i\mu_n)U_n\xrightarrow[~n\longrightarrow+\infty~]{}0~~\text{ in }X~. 
\end{equation}
Notice that, here, $u_n$ and $v_n$ are complex valued functions.

\vspace{3mm}

\noindent {\it $\bullet$ Second step: replacing $\gamma$ by a smooth damping.}\\
We recall that $H^1(\RR^d)$ is equipped with the convenient 
scalar product \eqref{prod-scalaire-H1}. Let us denote the
operator $\div(K(x)\grad)$ by $\Delta_K$. We have
$$
(A-i\mu_nId)U_n=\left(\begin{array}{c} v_n-i\mu_nu_n\\
(\Delta_K-Id)u_n  - \gamma(x)v_n-i\mu_nv_n \end{array}\right)
$$
and 
$$\Re(\langle (A-i\mu_n)U_n|U_n\rangle_X)=-\int \gamma(x)|v_n(x)|^2\,dx~.$$
Thus, \eqref{eq-absurde} implies that $\int \gamma(x)|v_n(x)|^2\,dx$ goes to
zero. Therefore, we can replace $\gamma$ by  any smooth damping
$\underline{\gamma}$ satisfying $0\leq \underline{\gamma}
\leq \gamma$ without changing \eqref{eq-absurde}. 
Let us show that we can construct such a damping
$\underline{\gamma}\in\Cc^\infty_b(\RR^d)$. First
choose $\theta=\max(0,\gamma-\varepsilon)$. For a small enough
$\varepsilon>0$, the damping $\theta$ still satisfies that its mean value
$\langle \theta \rangle_T(x,\xi)$ is uniformly bounded away from $0$. Moreover,
since $\gamma$ is uniformly continuous, the support of $\theta$ stays
at a uniform distance $\delta>0$ of the set where $\gamma$
vanishes. Now, mollify $\theta$ into $\theta * \rho_{\delta}$ where
$\rho_\delta$ is a $\Cc^\infty$ regularisation kernel with support in
$B(0,\delta)$. We obtain a smooth damping $\underline{\gamma}$ with a support
included in the one of $\gamma$. Thus, one can use this new damping
without changing \eqref{eq-absurde}. Moreover, this new damping
$\underline{\gamma}$ belongs to $\Cc^\infty_b(\RR^d)$, which ensures that the
multiplication by $\underline{\gamma}$ is a pseudo-differential operator of
order $0$. In the remaining part of this proof, to simplify the notations,
we will assume that $\gamma$ itself belongs to $\Cc^\infty_b(\RR^d)$.

\vspace{3mm}

\noindent {\it $\bullet$ Third step: separation between high and low
frequencies.}\\
We now work with a smooth damping $\gamma$ with bounded
derivatives. To obtain a contradiction from \eqref{eq-absurde}, we consider two
cases.
\begin{itemize}
\item[$\rightarrow$] {\bf High frequencies:} assume that $|\mu_n|$ goes to
$+\infty$. Since $A$ is a real operator, by symmetry, we can assume that
$\mu_n>0$ and we set $h_n=1/\mu_n$. We have to show that one cannot have 
$\|U_n\|_X=1$ and $(A-i/h_n)U_n\longrightarrow 0$. This will be shown in
Section \ref{sect-high} by semiclassical pseudo-differential arguments, using
the geometric control condition of Theorem \ref{th-main}.
\item[$\rightarrow$] {\bf Low frequencies:} assume that $(\mu_n)$ has a bounded
subsequence. Then, up to extracting a subsequence, one can assume that
$(\mu_n)$ converges to a real number $\mu$. Then \eqref{eq-absurde} is
equivalent to have a sequence $(U_n)$ with $\|U_n\|_X=1$ and
$(A-i\mu)U_n\longrightarrow 0$. In Section \ref{sect-low}, we will show that
this is not possible by using a global Carleman estimate. In this part, it is
in fact sufficient to replace the geometric control condition by the 
network control condition (NCC) stated in Theorem \ref{Th.2}. A similar
argument was developped independently by Le Rousseau and Moyano for the study of
the Kolmogorov equation.
\end{itemize}
Since Sections \ref{sect-high} and \ref{sect-low} provide a contradiction in
both cases, Theorem \ref{th-Huang} yields the proof of Theorem \ref{th-main}.

%%%%%%%%%%%%%%%%%%%%%%%%%%%%%%%%%%%%%%%%%%%%
%%%%%%%%%%%%%%%%%%%%%%%%%%%%%%%%%%%%%%%%%%%%

\section{Proof of Theorem~\ref{th-main}: high frequencies}\label{sect-high}

The purpose of this section is to obtain a contradiction from the existence of
sequences $(U_n)$ with $\|U_n\|_X=1$ and $(h_n)$ with $h_n\rightarrow 0$
satisfying $(A-i/h_n)U_n\longrightarrow 0$. To simplify the notations, we may
forget the index $n$ for the remaining part of this section and set
$U_n=U_h=(u_h,v_h)$. We have
$$
\left\{\begin{array}{l}v_h-\frac ih u_h=o_{H^1}(1)\\
(\Delta_K-id)u_h-\gamma(x)v_h-\frac ih v_h=o_{L^2}(1) \end{array}\right. 
$$
and thus
\begin{equation}\label{hf1}
\left\{\begin{array}{l}v_h-\frac ih u_h=o_{H^1}(1)\\
h^2(\Delta_K-Id)u_h-ih\gamma(x)u_h+u_h=o_{L^2}(h^2)+o_{H^1}(h)
\end{array}\right. 
\end{equation}
To obtain a contradiction between
\eqref{hf1} and Hypothesis (GCC) of Theorem \ref{th-main}, we will use
the semiclassical microlocal analysis and follow the ideas of the chapter 5 of
\cite{Zworski}. Notice that the usual way to deal with high frequencies is to
use semiclassical defect measures (see for example \cite{Zworski}). However,
this is not possible in our case since we work in an unbounded
domain and the semiclassical defect measure will only tell us what happens in
compact subsets.

\begin{lemma}\label{lemme-1-hf}
Assume that the operator 
$P_h=h^2(\Delta_K-Id)-ih\gamma(x)+Id$
has a resolvent in $L^2(\RR^d)$ satisfying
\begin{equation}\label{aa}
\|(P_h)^{-1}f\|_{L^2}\leq \frac Ch \|f\|_{L^2}~.
\end{equation}
Then $P_h$ has a resolvent in $H^1(\RR^d)$ also satisfying
$$\|(P_h)^{-1}f\|_{H^1}\leq \frac Ch \|f\|_{H^1}~.$$
\end{lemma}
\begin{demo}
We argue by contradiction. Assume that there exists a sequence $(u_h)$ with
$\|u_h\|_{H^1}=1$ and $P_hu_h=o_{H^1}(h)$. Multiplying by $\overline{u_h}$ and
integrating, we get that 
$$-h^2\|u_h\|_{H^1}^2-ih\int_{\RR^d}\gamma(x)|u_h|^2+\|u_h\|_{L^2}^2
=o(h)\|u_h\|_{L^2}~.$$
Taking the real part and solving the equation in $\|u_h\|_{L^2}$, we get
$$\|u_h\|_{L^2}=\frac12\left(o(h)+\sqrt{o(h^2)+4h^2\|u_h\|_{H^1}^2}\right)\sim
h~.$$
We introduce the operator $\grad_K=(-\Delta_K+Id)^{1/2}$. It is the particular
case $h=1$ of the semiclassical operator $(-h^2\Delta_K+Id)^{1/2}$, which has
for principal symbol $\sqrt{\xi\t.K(x).\xi+1}$ (see
Section \ref{sect-A-pseudo} in Appendix for a brief recall about
pseudo-differential semiclassical calculus). Obviously, it commutes with any
polynomial of $\Delta_K$. Moreover, applying i) of Corollary \ref{coro-commut}
of the Appendix, with $h=1$ fixed, we get that the commutator
$[\grad_K,\gamma(x)\cdot]$ is bounded in $L^2(\RR^d)$. Thus, since
$u_h=\mathcal{O}_{L^2}(h)$,
$$\grad_K (P_hu_h)=P_h(\grad_K u_h)-ih[\grad_K,\gamma(x)]u_h=P_h(\grad_K
u_h)+\mathcal{O}_{L^2}(h^2)~.$$
Since $P_hu_h=o_{H^1}(h)$, we obtain that $P_h(\grad_K u_h)=o_{L^2}(h)$. Using
the assumption on the resolvent of $P_h$, we obtain that $\grad_K u_h$ goes to
$0$ in $L^2(\RR^d)$ when $h$ goes to $0$. However, $\|\grad_K u_h\|_{L^2}$
is equivalent to $\|u_h\|_{H^1}$ and we obtain a contradiction with the
assumption $\|u_h\|_{H^1}=1$.
\end{demo}

\begin{prop}\label{prop-1-hf}
If the operator $P_h=h^2(\Delta_K-Id)-ih\gamma(x)+Id$
has a resolvent in $L^2(\RR^d)$ satisfying \eqref{aa}, then \eqref{hf1} cannot
hold.
\end{prop}
\begin{demo}
We argue by contradiction again. Assume that $P_h$ satisfies \eqref{aa} and
assume that there exists $U_h=(u_h,v_h)$ with $\|U_h\|_X=1$ such that
\eqref{hf1} holds. As in the beginning of the proof of Lemma
\ref{lemme-1-hf}, multiplying the second equation of \eqref{hf1} by
$\overline u_h$, integrating, taking the real part and solving the equation of
second degree in $\|u_h\|_{L^2}$, we get that
$$\|u_h\|_{L^2}=\frac12\left(o(h)+\sqrt{o(h^2)+4h^2\|u_h\|_{H^1}^2}\right)~.$$
Due to the first equation of \eqref{hf1} and since $\|U_h\|_X=1$, we must have
$\|u_h\|_{H^1}\sim 1/\sqrt{2}$, $\|u_h\|_{L^2}\sim h/\sqrt{2}$ and
$\|v_h\|_{L^2}\sim 1/\sqrt{2}$. 

We introduce $w_h=P_h^{-1}(f_h)$ where $f_h$ is the term $o_{H^1}(h)$ in the
second equation of \eqref{hf1}. By assumption and by Lemma \ref{lemme-1-hf}, we
have that $w_h=o_{H^1}(1)$. By the same straightforward computation than the
one just above, we also have $w_h=o_{L^2}(h)$. Then, $u_h-w_h$ solves
$P_h(u_h-w_h)=o_{L^2}(h^2)$ and by assumption, we get that $u_h-w_h=o_{L^2}(h)$
and thus that $u_h=o_{L^2}(h)$. This is a contradiction with the fact that
$\|u_h\|_{L^2}\sim h/\sqrt{2}$, which was proved above.
\end{demo}

Due to Proposition \ref{prop-1-hf}, to obtain a contradiction from \eqref{hf1},
it remains to show the $L^2$-resolvent estimate \eqref{aa}. Obtaining this
estimate is the central argument for controlling the high frequencies. Here, we
will use pseudo-differential calculus and we will see the importance of
Hypothesis (GCC) of Theorem \ref{th-main}. The remaining part of this section is
thus devoted to the proof of the following result.
\begin{prop}\label{prop-hf}
The operator 
$$P_h=h^2(\Delta_K-Id)-ih\gamma(x)+Id$$
has a resolvent in $L^2(\RR^d)$ satisfying
$$\|(P_h)^{-1}f\|_{L^2}\leq \frac Ch \|f\|_{L^2}~.$$
\end{prop}
\begin{demo}
As usual, we argue by contradiction and assume that there exists a sequence
$(h_n)$ going to zero and functions $(u_n)\subset H^2(\RR^d)$ such that
$\|u_n\|_{L^2}=1$ and $P_{h_n}u_n=o_{L^2}(h_n)$. Once again, we may forget the
indices and assume that $\|u_h\|_{L^2}=1$ and 
\begin{equation}\label{hf2}
h^2(\Delta_K-Id)u_h-ih\gamma(x)u_h+u_h=o_{L^2}(h)~.
\end{equation}
In what follows, we will use the notations and the results of the
pseudo-differential semiclassical calculus recalled in Section
\ref{sect-A-pseudo}. Our proof follows the lines of Chapter 5 of
\cite{Zworski}, omitting the notion of defect measures, which is not convenient
in the case of unbounded domains.

\vspace{3mm}

\noindent {\it $\bullet$ First step: $u_h$ is concentrating along the radial
speeds $\xi\t K(x)\xi=1/h^2$.}\\
First notice that the main part of $P_h$ is $h^2\Delta_K+Id$ in the sense
that $(h^2\Delta_K+Id)u_h=o_{L^2}(1)$. As explained in Appendix, up to an
error term $\mathcal{O}_{L^2}(h^2)$, this main part is a pseudo-differential
semiclassical operator $\Op_h(-\xi\t K(x)\xi+1)$.
Let $\chi(x,\xi)\in\Cc^\infty(\RR^d,\RR_+)$ be a smooth cutting function which
is equal to one in a neighbourhood of the sphere $\Sigma=\{(x,\xi),\xi\t
K(x)\xi=1\}$ and equal to $0$ outside the annulus  $1/2K_{\rm sup}\leq
|\xi| \leq 2/K_{\rm inf}$. Also assume that $\chi$ and its derivatives are
bounded,
which implies that $\chi(x,\xi)$ is a symbol of order $0$. We claim that $u_h$
is concentrating on the microlocal set $\{(x,\xi),\xi\t K(x)\xi=1/h^2\}$ in the
sense that
$\langle \Op_h(1-\chi(x,\xi))u_h|u_h\rangle_{L^2}$ goes to $0$ when $h$ goes to
$0$.

To prove this claim, we introduce another
smooth cutting function $\theta$ which is equal to $1$ in a neighbourhood of the
sphere $\Sigma=\{(x,\xi),\xi\t K(x)\xi=1\}$ and equal to $0$ in the support of
$1-\chi$. The symbol
$a(x,\xi)=-\xi\t K(x)\xi+1+i\theta$ is of order $2$ and uniformly bounded away
from $0$. By Corollary \ref{coro-inverse} in Appendix, the symbol 
$b(x,\xi)=\frac1{a(x,\xi)}$ is of order $-2$ and satisfies 
$$\Op_h(a)\circ\Op_h(b)=Id+\mathcal{O}_{L^2}(h)~~~\text{ and }~~~
\Op_h(b)\circ\Op_h(a)=Id+\mathcal{O}_{L^2}(h)~.$$
Thus, 
$$\langle \Op_h(1-\chi)u_h|u_h\rangle_{L^2}=\langle
\Op_h(1-\chi)\circ\Op_h(b)\circ\Op_h(a)u_h|u_h\rangle_{L^2}+\mathcal{O}(h)~.$$
On the other hand, $\Op_h(a)=\Op_h(-\xi\t K(x)\xi+1)+i\Op_h(\theta)$ and thus 
$\Op_h(a)u_h=o_{L^2}(1)+i\Op_h(\theta)u_h$. Since $1-\chi$ and $b$ are of order
$0$ or less, their corresponding operators are bounded in $L^2(\RR^d)$,
uniformly with respect to $h$ and 
$$\langle
\Op_h(1-\chi)u_h|u_h\rangle_{L^2}=i\langle\Op_h(1-\chi)\Op_h(b)\Op_h(\theta)
u_h|u_h\rangle_{L^2}+o(1)~.$$
Now, it remains to apply Proposition \ref{prop-compo} in Appendix to see
that, since $1-\chi$ and $\theta$ have disjoint supports,
$$\Op_h(1-\chi)\Op_h(b)\Op_h(\theta)u_h=\Op_h((1-\chi)b\theta)+o_{L^2}(1)=o_{L^2
} (1)~.$$
This shows that 
$$\langle \Op_h(1-\chi)u_h|u_h\rangle_{L^2} \xrightarrow[~h\longrightarrow
0~]{}0~.$$

\vspace{3mm}

\noindent {\it $\bullet$ Second step: using the geometric control condition of
Theorem \ref{th-main}.}\\
First notice that 
$$P_h=\Op_h(-\xi\t K(x)\xi+1)-ih\Op_h(\gamma(x))+\mathcal{O}_{L^2\rightarrow
L^2} (h^2)$$ 
and that we may assume that $\gamma$ is smooth and bounded and so that it is a
symbol of order $0$ (see Section \ref{sect-proof}). Let $a(x,\xi)$ be a
symbol of order $0$. By Corollary \ref{coro-commut} in Appendix,
the commutator of $P_h$ and $\Op_h(a)$ is 
$$[\Op_h(a),P_h]=-ih\Op_h\big(\{\xi\t K(x)\xi,a(x,\xi)\}\big)+\mathcal{O}_{L^2
\rightarrow L^2} (h^2)~.$$
On the other hand, since $P_hu_h=o_{L^2}(h)$,
\begin{align*}
\langle [\Op_h(a),P_h]u_h|u_h\rangle_{L^2}&=\langle
\Op_h(a)P_hu_h|u_h\rangle_{L^2}-\langle
P_h\Op_h(a)u_h|u_h\rangle_{L^2}\\
&=o(h)-\langle \Op_h(a)u_h|P_h^*u_h\rangle_{L^2}\\
&=-\langle \Op_h(a)u_h|(P_h+2ih\gamma(x))u_h\rangle_{L^2}+o(h)\\
&=2ih \langle\Op_h(a)u_h|\gamma(x)u_h\rangle_{L^2}+o(h)\\
&=2ih \langle\gamma(x)\Op_h(a)u_h|u_h\rangle_{L^2}+o(h)\\
&=2ih \langle\Op_h(a\gamma)u_h|u_h\rangle_{L^2}+o(h)
\end{align*}
Thus, setting $g(x,\xi)=\xi\t K(x)\xi$, we obtain that
\begin{equation}\label{hf3}
\langle\Op_h(2a\gamma+\{g,a\})u_h|u_h\rangle_{L^2}\xrightarrow[~h\longrightarrow
0~]{}0~.
\end{equation}
Due to Corollary \ref{coro-inverse} of Appendix, we will get a
contradiction
with $\|u_h\|_{L^2}=1$ if we find $a$ such that $2a\gamma+\{g,a\}$ is uniformly
bounded away from zero. Assume that $a(x,\xi)$ is constant equal to $1$ for
large $\xi$, then $2a\gamma+\{g,a\}$ is a symbol of order $0$. Moreover, the
first step of this proof shows that modifying $2a\gamma+\{g,a\}$ away from the
sphere $\Sigma=\{(x,\xi),\xi\t K(x)\xi=1\}$ has no influence on \eqref{hf3}.
Thus, it is sufficient
to exhibit a symbol $a$  such that $2a\gamma+\{g,a\}$ is uniformly
bounded and stay uniformly away from zero on $\Sigma$. 

Let us recall that $\varphi_t$ is the Hamiltonian flow associated to $g$ and
that $T$ is a time such that the mean value $\langle \gamma
\rangle_T(x,\xi)=\frac1T\int_0^T \gamma(\varphi_t(x,\xi))dt$ is uniformly
bounded away from $0$ away from $\Sigma$, according to Assumption (GCC) of
Theorem \ref{th-main}. We choose $a(x,\xi)=e^{c(x,\xi)}$ with 
$$c(x,\xi)=\frac2T \int_0^T(T-t)\gamma(\varphi_t(x,\xi))\,dt=\frac2T
\int_0^T\int_0^t\gamma(\varphi_s(x,\xi))\,ds\,dt~.$$
By definition of the Hamiltonian flow, for any function
$f\in\Cc^1(\RR^{2d},\RR)$, we have 
$$\{g,f\}(x,\xi)=\partial_\tau f(\varphi_\tau(x,\xi))_{|\tau=0}~.$$
Since
\begin{align*}
c(\varphi_\tau(x,\xi))&=\frac2T
\int_0^T(T-t)\gamma(\varphi_{t+\tau}(x,\xi))\,dt\\
&=\frac2T \int_\tau^{T+\tau}(T-t+\tau)\gamma(\varphi_t(x,\xi))\,dt
\end{align*}
we get that
$$\{g,c\}(x,\xi)=\frac 2T \int_0^T\gamma(\varphi_t(x,\xi))\,dt -2\gamma(x,\xi) =
2\langle \gamma\rangle_T(x,\xi)-2\gamma(x,\xi)~.$$
Thus, we have
$$2a\gamma+\{g,a\}=2e^{c(x,\xi)}\langle
\gamma\rangle_T(x,\xi)~.$$
By assumption (GCC) of Theorem \ref{th-main} and since $c\geq 0$, there exists
$\alpha>0$ such that, for all $(x,\xi)\in\Sigma$,
$2a\gamma+\{g,a\}\geq \alpha>0$. As explained above, we can neglect any
$(x,\xi)$ away from $\Sigma$ and this yields that
$$\langle\Op_h(2a\gamma+\{g,a\})u_h|u_h\rangle_{L^2} \sim
\langle\Op_h\left(2e^{c(x,\xi)}\langle
\gamma\rangle_T(x,\xi)\right)u_h|u_h\rangle_{L^2} \geq 2\alpha \|u_h\|^2_{L^2}~,
$$
which contradicts \eqref{hf3} since $\|u_h\|_{L^2}=1$.
\end{demo}

%%%%%%%%%%%%%%%%%%%%%%%%%%%%%%%%%%%%%%%%%%%%
%%%%%%%%%%%%%%%%%%%%%%%%%%%%%%%%%%%%%%%%%%%%

\section{Proof of Theorem~\ref{th-main}: low frequencies}\label{sect-low}
We now have to deal with the low frequencies to finish the proof of
Theorem~\ref{th-main}. This is done by using Carleman estimates. Notice that
the same tool will provide the logarithmic decay of Theorem \ref{Th.2} (see
Section \ref{S-helmoltz}).

In this section, we fix first  a real number $\mu$ and we 
assume that there is a sequence $(U_n)$ with $\|U_n\|_X=1$ and
$(A-i\mu)U_n\longrightarrow 0$, that is that $U_n=(u_n,v_n)$ satisfies
$v_n=i\mu u_n+o_{H^1}(1)$ and 
\begin{equation}\label{lf1}
(\Delta_K-Id)u_n-i\mu\gamma(x)u_n+\mu^2 u_n=o_{L^2}(1)~.
\end{equation}
We work with $\gamma\in\Cc^\infty_b(\RR^d)$ satisfying the geometric
control condition (GCC) of Theorem \ref{th-main} (see Section \ref{sect-proof}).
This condition yields that, for any $(x,\xi)\in\Sigma$, the ray of length $T$
contains a point $y$ such that $\gamma (y)\geq \alpha/T$. Since the metric $K$
is uniformly bounded, we can find a sequence $(x_n)\subset \RR^d$ such that
$\gamma (x_n)\geq \alpha/T$ and
any point of $\RR^d$ is at bounded distance $L$ of the set $\cup_n \{x_n\}$.
Since $\gamma$ is uniformly continuous, we can find $r,a>0$ such that
$\gamma(x)\geq a >0$ on $\cup_n B(x_n,r)$. This is the control by a network of
balls (NCC) stated in Theorem \ref{Th.2}, which will be a sufficient
condition to control the low frequencies in this section. We will denote by
$\omega$ the set $\omega=\cup_n B(x_n,r)$.

As a first basic computation, we can multiply \eqref{lf1} by $\overline u_n$ and
integrate. Taking real and imaginary parts,  we obtain
\begin{equation}\label{lf2}
\|u_n\|^2_{H^1}=\mu^2 \|u_n\|^2_{L^2}+o(1)~~~~~~\text{ and }~~~~~~
\int_{\RR^d}\gamma(x)|u_n(x)|^2\,dx \xrightarrow[~n\longrightarrow \infty~]{}
0~. 
\end{equation}
Also notice that $\|v_n\|_{L^2}=\mu \|u_n\|_{L^2}+o(1)$. Thus, if $\mu=0$ we
obtain a contradiction between $\|U_n\|_X=1$ and $\|U_n\|_X \sim
\mu\|u_n\|_{L^2} + o(1)$. Assume from now on that $\mu\neq 0$, we get that
$\|U_n\|_X$ is equivalent to $\|u_n\|_{L^2}$ and, up to a renormalisation, we
can assume that $(u_n)$ satisfies $\|u_n\|_{L^2}=1$. 

Now, we can multiply \eqref{lf1} by $\Delta_K \overline u_n$ and
integrate. The real part of the result shows that $\|\Delta_K
u_n\|^2=\Oc(1)+\langle o_{L^2}(1) | \Delta_K u_n\rangle_{L^2}$, which
implies that $\|\Delta_K u_n\|$ is bounded. Then, considering the imaginary
part, we get that  
\begin{equation}\label{lf3}
\int_{\RR^d}\gamma(x)|\grad u_n(x)|^2\,dx \xrightarrow[~n\longrightarrow
\infty~]{} 0~. 
\end{equation}

\vspace{3mm}

\noindent {\it $\bullet$ First step: using H\"ormander sub-ellipticity
argument.}\\
We set $P=\Delta_K-i\mu\gamma(x)+(\mu^2-1)Id$ and 
$$Q_h=h^2 e^{\varphi/h} (-\Delta_K+(1-\mu^2)Id) e^{-\varphi/h}~.$$
We follow the classical arguments (see for example \cite{Lerousseau}). We have
$$Q_h u=-h^2\Delta_K u+2h\grad\varphi\t K(x) \grad u + u(
 - \grad\varphi\t K(x)\grad \varphi+ h\Delta_K(\varphi) + h^2(1-\mu^2))~.$$
With the notations of the Appendix \ref{sect-A-pseudo}, using
$-h^2\Delta_K=\Op_h( \xi\t K(x)\xi )$, $h\grad = \Op_h(i\xi)$ and Proposition
\ref{prop-compo}, we obtain that 
$$Q_h = \Op_h\left( \xi\t K(x) \xi - \grad \varphi\t K(x) \grad\varphi + 2 i
\grad \varphi\t K(x)\xi\right) + \Oc_{L^2\rightarrow L^2}(h^2)~.$$
We set $Q_h= Q^R_h+i Q^I_h +  \Oc_{L^2\rightarrow L^2}(h^2)$ with 
\begin{align*}
Q^R_h&=\Op_h(q_R)=\Op_h (  \xi\t K(x) \xi - \grad \varphi\t K(x) \grad\varphi )
\\
Q^I_h&=\Op_h(q_I)=\Op_h ( 2 \grad \varphi\t K(x)\xi )~.
\end{align*}

We use Proposition \ref{prop-adjoint} in Appendix to check that $Q^R_h$ and
$Q^I_h$ are self-adjoint operators and Corollary \ref{coro-commut} shows that  
\begin{align}
 \|Q_hu\|^2_{L^2}&=  \|Q^R_h u \|^2_{L^2} + \|Q^I_h u \|^2_{L^2} + \langle
Q^R_h u | iQ^I_h u\rangle_{L^2} + \langle iQ^I_h u | Q^R_h u\rangle_{L^2} +
\Oc(h^2\|u\|^2_{L^2}) \nonumber\\
&=  \langle ((Q^R_h)^2 + (Q^I_h)^2 + i[Q^R_h,Q^I_h]) u | u\rangle_{L^2}
+ \Oc(h^2\|u\|^2_{L^2})\nonumber\\
&\geq h \langle \Op_h ( \eta (q_R^2+q_I^2)+\{q_R,q_I\}) u | u \rangle +
\Oc(h^2\|u\|^2_{L^2}) \label{eq-demo-123}
\end{align}
where $\eta$ is any number such that $\eta h\leq 1$. Let
$$B=\left\{\,(x,\xi)\in\RR^{2d}~,~~  \frac{K_{\rm inf}}{2K_{\rm sup}} |\grad
\varphi(x)| \leq  |\xi| \leq  \frac{2K_{\rm sup}}{K_{\rm inf}} |\grad
\varphi(x)|\,\right\}~.$$
Notice that $q_R$ is uniformly away from $0$ for $(x,\xi)$ outside $B$.
Assume that $\varphi$ satisfies the sub-ellipticity criterion:
There exists $\alpha>0$ such that 
\begin{equation}\label{eq-Hormander}
\{q_R,q_I\}(x,\xi)\,\geq\,\alpha\,>\,0\text{ on }((\RR^d\setminus\omega)\times
\RR^d) \cap \{(x, \xi);  q_R (x, \xi) = q_I(x, \xi) =0\}  
\end{equation}
Then, taking $\eta$ sufficiently large, $\eta (q_R^2+q_I^2)+\{q_R,q_I\}$ is
uniformly positive on $\RR^d\setminus\omega$. Moreover, the behaviour for large
$|\xi|$ is given by $q_R^2$, thus we have that $\eta (q_R^2+q_I^2)+\{q_R,q_I\}$
is a symbol of order $4$ and there is a positive constant $\alpha>0$ such that 
$$\eta (q_R^2+q_I^2)+\{q_R,q_I\} \geq \alpha (1+|\xi|^2)^2  ~~\text{ on
}(\RR^d\setminus \omega)\times\RR^d~.$$
Using G\aa{}rding inequality stated in Proposition \ref{prop-garding-2} in
Appendix and \eqref{eq-demo-123}, we obtain that, if \eqref{eq-Hormander} holds,
then there is $c>0$ such that 
\begin{equation}\label{eq-demo-124}
\|Q_hu\|^2_{L^2}~\geq~ch\|u\|^2_{L^2}~~~~\text{ for all }u\text{ satisfying
}u_{|\omega}\equiv 0~.
\end{equation}
To obtain a contradiction with \eqref{lf1} and $\|u_n\|_{L^2}=1$, we proceed
as follows. Let $\chi\in\Cc^\infty_b(\RR,[0,1])$ be a function such that
$\chi(s)=0$ for $s\geq a$ and $\chi\equiv 1$ in a neighbourhood of $0$. We have
that $\chi\circ \gamma$ vanishes on $\omega$ and, since $(1-\chi\circ\gamma)$
vanishes on $\{x,\gamma(x)\leq \nu\}$ for some small $\nu>0$, any derivative of
$\chi\circ\gamma$ is controlled by $\kappa\gamma$ for $\kappa$ large enough. We
set
$v_n=e^{\varphi/h}(\chi\circ\gamma)u_n$ and compute 
\begin{align*}
Q_hv_n~&=~ h^2 e^{\varphi/h} (-P-i\mu\gamma)(\chi\circ\gamma)u_n \\
&=~ -h^2e^{\varphi/h}(\chi\circ\gamma)Pu_n ~-~ i\mu
h^2e^{\varphi/h}\gamma(\chi\circ\gamma)u_n~~-2h^2e^{\varphi/h}
\grad(\chi\circ\gamma)\t K(x) \grad u_n\\
&~~~~~-~h^2e^{\varphi/h} u_n \Delta_K
(\chi\circ\gamma)\\
&=~o_{L^2}(1)~~~~~~~~~~~~\text{ when }n\rightarrow +\infty\text{ and }h>0
\text{ is fixed.}
\end{align*}
where we used the fact that $\gamma u_n$ and $\gamma \grad u_n$ goes to zero in
$L^2(\RR^d)$. Now, remember that $u_n - (\chi\circ \gamma)u_n$ is supported on
$\{x,\gamma(x)\geq \nu\}$ and thus also goes to
zero in $L^2$. For $h$ fixed, $\|v_n\|_{L^2}$ is thus uniformly positive since
$e^{\varphi/h}(x)\geq e^{\min \varphi/h}>0$ and
$\|(\chi\circ\gamma)u_n\|_{L^2}\rightarrow 1$. Since $v_n$ vanishes on $\omega$,
this is an obvious contradiction with \eqref{eq-demo-124} and
$Q_hv_n=o_{L^2}(1)$ shown above.

\vspace{3mm}

\noindent {\it $\bullet$ Second step: Carleman weight $\varphi=e^{\lambda
\psi}$.}\\
The usual way to obtain a weight $\varphi$ satisfying
H\"ormander sub-elliptic assumption \eqref{eq-Hormander} consists in
choosing $\varphi=e^{\lambda \psi}$, with a function $\psi$, whose critical
points belongs to $\omega=\cup_n B(x_n,r)$, and $\lambda$ large enough. We
reproduce here this argument with an obvious care about uniformity.

Assume that $\varphi=e^{\lambda\psi}$ for some constant
$\lambda$ and that $\psi\in\Cc^\infty_b(\RR^d)$ is such that there exists
$\alpha>0$ such that $|\grad \psi(x)|\geq \alpha >0$
for all $x\in\RR^d\setminus \omega$. A straightforward computation yields
\begin{align*}
\{q_R,q_I\} & = \partial_\xi\left(\xi\t K(x) \xi\right)
\partial_x\left(\lambda e^{\lambda\psi}\grad\psi \t K(x)\xi\right)\\
&~~~~~~~~- \partial_x\left(\xi\t K(x)
\xi - \lambda^2 e^{2\lambda \psi} \grad\psi\t K(x)\grad
\psi\right)\partial_\xi\left(\lambda e^{\lambda\psi}\grad\psi \t K(x)\xi
\right)\\
&\geq \lambda^4 e^{3\lambda \psi}(\grad\psi\t K(x) \grad\psi)^2 + \Oc
(|\xi|^2\lambda e^{\lambda\psi}) + \Oc(\lambda^3 e^{3\lambda\psi})
\end{align*}
where the estimations $\Oc(\cdot)$ hold for $|\xi|$ and $\lambda$ going to
$+\infty$. Notice that, when $(x,\xi)$ belongs to $B$, $|\xi|$ is of order
$\Oc(\lambda e^{\lambda \psi(x)})$. Since $|\grad
\psi(x)|\geq \alpha >0$ on $\RR^d\setminus \omega$, we can fix $\lambda$ large
enough, such that the positive term $\lambda^4 e^{3\lambda \psi}(\grad\psi\t
K(x) \grad\psi)^2$ controls the last two terms (with indefinite sign), when
$(x,\xi)$ belongs to $B$ and $x\not\in\omega$.

\vspace{3mm}

\noindent {\it $\bullet$ Third step: construction of the appropriate Carleman
phase $\psi$.}\\
To summarise, the above arguments show that, if we are able to
construct a suitable phase $\psi$, then the sub-elliptic criterium 
\eqref{eq-Hormander} would hold and also the uniform positivity property
\eqref{eq-demo-124}. This will provide a contradiction between \eqref{lf1} and
$\|u_n\|_{L^2}\equiv 1$, as described below Equation \eqref{eq-demo-124}. This
will yield the control of the low frequencies and finished this section. 

Thus, we only have to construct $\psi\in\Cc^\infty_b(\RR^d,\RR)$ such that
 $|\grad \psi|$ is uniformly positive outside $\omega=\cup_n
B(x_n,r)$. We recall that the points $x_n$ are assumed to form a network in the
sense that any point $x\in\RR^d$ is at distance at most $L$ of a point $x_n$. We
split $\RR^d$ in cubes $(C_k)_{k\in\ZZ^d}$ of size $4L$ by setting
$C_k=4Lk+[-2L,2L]^d$.
For each $k\in\ZZ^d$, the center of $C_k$ is $c_k=4Lk$ and there is at least
one of the points $x_n$ which is in $c_k+[-L,L]^d$, let us denote it $y_k$. For
each $k$, one can find a $\Cc^\infty-$diffeomorphism with compact support in
the interior of $C_k$ such that $y_k$ is mapped onto $c_k$. We glue all these
diffeomorphisms into a diffeomorphism $\Phi$ of $\RR^d$ mapping all the $y_k$
onto the $c_k$ and we notice that we can make this construction such that $\Phi$
and $\Phi^{-1}$ belong to $\Cc^\infty_b(\RR^d,\RR)$ (an explicit construction
is given in \cite{LRM}). Consider 
$$\psi(x)=\tilde\psi\circ \Phi(x)~~~~\text{ with }~~~
\tilde\psi(x)=\sum_{i=1}^d
\cos\left(\frac{\pi x_i}{4L}\right)~.$$
Obviously, $|\grad \tilde\psi|$ is uniformly positive outside $\cup_n 
B(c_n,\rho)$ for any $\rho>0$. Thus, $|\grad \psi|$ is uniformly positive
outside $\cup_n B(x_n,r)$, which concludes this section.

%%%%%%%%%%%%%%%%%%%%%%%%%%%%%%%%%%%%%%%%%%%%
%%%%%%%%%%%%%%%%%%%%%%%%%%%%%%%%%%%%%%%%%%%%

\section{Proof of Theorem~\ref{Th.2}}\label{S-helmoltz}

The proof of Theorem \ref{Th.2} relies on the same arguments than the
ones of Sections \ref{sect-proof} and \ref{sect-low}. We will only outline
them. 

Instead of Theorem \ref{th-Huang}, we use the following characterisation of the
logarithmic rate of decay given by Theorem~3 of Burq~\cite{Burq}.
\begin{theorem}[Burq, \protect{\cite[Theorem 3]{Burq}}] \label{th-burq}
Let  $A$ be maximal dissipative operator (and hence the generator of $e^{At}$
a $\Cc^0-$semigroup of contractions) in a Hilbert space $X$ and assume that
there exist $C,c>0$ such that $i\RR\subset \rho(A)$ and 
\begin{equation}\label{eq-burq}
\forall \mu \in \mathbb{R},  \|(A-i\mu Id)^{-1}\|_{\Lc(X)}~<C e^{c|\mu|}.
\end{equation}
Then for any $k>0$ there exists $C_k$ such that for any $t>0$,
$$ \left\| \frac{ e^{tA}} { (1- A)^k} \right\| _{\mathcal{L} (X)} \leq
\frac{C_k} { \log(2+t)^k}.$$
\end{theorem}

First notice that the estimate \eqref{eq-burq} is already proved for low
frequencies in Section \ref{sect-low}. Thus, it is enough to prove it for large
$\mu$ (say $\mu = h^{-1}, h \rightarrow 0^+$, the case of negative $\mu$ being
similar). We are consequently in a high frequency regime, but are nevertheless
going to use the approach developed in the previous section  for low
frequencies, based on Carleman weight and H\"ormander sub-ellipticity
argument. Indeed, a simple adaptation allows to prove similar Carleman
estimates in the high-frequency regime (and hence for the semi-classical
Helmoltz operator) by tracking the exponential dependence of the constants
with respect to the frequency parameter.

As in Section \ref{sect-proof}, we can assume that $\gamma$ is smooth by
arguing by contradiction: assume that there
exist two sequences $(U_n)=(u_n,v_n)\subset D(A)=H^2(\RR^d)\times H^1(\RR^d)$
and $(\mu_n)\longrightarrow +\infty$ such that $\|U_n\|_X^2=1$ and
$\|(A-i\mu_n)U_n\|_X$ goes to $0$ in $X$ faster than any exponential. Once
again, we have $\Re(\langle (A-i\mu_n)U_n|U_n\rangle_X)=-\int
\gamma(x)|v_n(x)|^2\,dx$, which shows that $\gamma v_n$ is decaying as fast as
$(A-i\mu_n)U_n$. As in the second step of Section \ref{sect-proof}, we can
replace $\gamma$ by a smooth damping $\underline{\gamma}\in\Cc^\infty_b(\RR^d)$
satisfying $\underline{\gamma}\leq \gamma$ without changing the fact that
$\|(A-i\mu_n)U_n\|_X$ goes to $0$ in $X$ faster than any exponential. Notice
that, if $\gamma$ satisfies Hypothesis (NCC) of Theorem \ref{Th.2}, then one can
easily construct a smooth damping $\underline{\gamma}\leq \gamma$ also
satisfying (NCC).

We would like to show \eqref{eq-burq} for large $\mu$. As previously, simple
calculations show that it is enough to prove a similar estimate on
$((\Delta_K-Id)-i\mu\gamma(x)+\mu^2)^{-1}$: 
\begin{equation}\label{est.resh}
\exists C, c>0~,~~ \forall \mu \in \mathbb{R}~,~~ 
|||((\Delta_K-Id)-i\mu\gamma(x)+\mu^2)^{-1}|||_{\Lc(L^2( \mathbb{R}^d))}\leq C
e^{c|\mu|}~.
\end{equation}

We use the arguments and the notations of Section \ref{sect-low}.

Let $(u,f)$ solutions to $((\Delta_K-Id)-i\mu\gamma(x)+\mu^2)u =f$, i.e.,
setting $h=1/\mu$, 
$$(h^2\Delta_K-ih\gamma(x)+(1-h^2))u= h^2 f~.$$

Let 
$$\widetilde{Q}_h= e^{\varphi/h} (-h^2\Delta_K+(h^2-1)Id) e^{-\varphi/h}~.$$
We have $\widetilde{Q}_h = \widetilde{Q}_h^R + i \widetilde{Q}_h^I + \Oc(h^2)$
with 
\begin{align*}
\widetilde{Q}^R_h&=\Op_h(\widetilde{q}_R)=\Op_h (  \xi\t K(x) \xi - \grad
\varphi\t K(x) \grad\varphi -1) \\
\widetilde{Q}^I_h&=\Op_h(\widetilde{q}_I)=\Op_h ( 2 \grad \varphi\t K(x)\xi )~.
\end{align*}
In this setting, we shall assume that the phase function $\varphi\in C^\infty
_b(\mathbb{R}^d)$ satisfies H\"ormander hypo-ellipticity condition: there exists
$\alpha>0$ such that 
 \begin{equation}
\{\widetilde{q}_R,\widetilde{q}_I\}(x,\xi)\,\geq\,\alpha\,>\,0\text{ on
}((\RR^d\setminus\omega)\times 
\RR^d)  \cap \{(x, \xi);  \widetilde{q}_R (x, \xi) = \widetilde{q}_I(x, \xi) =0\} 
\end{equation}
The same proof as in the previous section shows that, under this condition, if
$v$ vanished in $\omega$, then for $h>0$ small enough,  
$$ \| v \|_{L^2} \leq \frac C h \| \widetilde{Q}_h v\|_{L^2}.$$
Coming back to $u$, applying the previous estimate to $v = e^{\varphi (x) /h }
\tilde\chi u$ with a cutoff $\tilde\chi=\chi\circ\gamma$ as in the previous
section,
we get 
 $$\| u\| _{L^2} \leq Che^{c/h} \| \tilde\chi f - (2(\grad
\tilde\chi)\t K(x)\grad + \Delta_K \tilde \chi)u\|_{L^2} +
Ce^{c/h}\|\gamma \tilde\chi u \|_{L^2}~,$$ 
% $$ \| u\| _{L^2} \leq Ce^{c/h} \| \chi f + ([ \Delta_K, \chi\circ\gamma ]) u
% \|_{L^2}~,$$
where $c = \sup _{x, y \in \mathbb{R}^2} |\varphi (x) - \varphi (y)|$.
Remember that all the terms involving $u$ in the right-hand side are
controlled by $\gamma u$, itself being controlled by the usual computation $\int
\gamma |u|^2= -h \Re(\int fu)$. We can now proceed by contradiction and
conclude the proof of the estimate~\eqref{est.resh} exactly as in the previous
section: it is impossible to have sequences $(\mu_n)=(1/h_n)$, $(u_n)$ and
$(f_n)$ such that $\|u_n\|=1$, $h_n\rightarrow 0$ and $(f_n)$ goes to zero
faster than any exponential $e^{-\kappa \mu_n}$.

To conclude, it remains to construct the Carleman weight $\varphi=
e^\lambda \psi$ satisfying H\"ormander hypo-ellipticity condition. This is done
exactly as in Section \ref{sect-low}. 
Notice that the only difference here with the low frequency case is that
$\widetilde{q}_R = q_R -1$. However, during the construction, this
additional term $1$ generates terms which are of order $\lambda ^2 e^{\lambda
\Psi}$. Thus, this will not perturb the exponential bound of the estimate.

%%%%%%%%%%%%%%%%%%%%%%%%%%%%%%%%%%%%%%%%%%%%
%%%%%%%%%%%%%%%%%%%%%%%%%%%%%%%%%%%%%%%%%%%%

\section{Applications to other problems}\label{sect-applis}

The exponential decay of the linear semigroup $e^{At}$ is an essential
assumption for obtaining several dynamical properties of the damped wave
equations. In this section, we emphasise different results, which are
corollaries of Theorem \ref{th-main}. Each result was already
known with stronger assumptions implying the exponential decay of $e^{At}$.
Since we have obtained this decay with weaker conditions, we can improve on these
results.
\subsection{Stabilisation in uniformly local Sobolev spaces}
Theorem \ref{th-main} concerns the solutions of the damped wave equation with
finite energy. In an unbounded domain, a natural question would be  to consider
solutions with infinite energy. For this reason, we introduce the uniformly
local Sobolev spaces as follows. For any $u \in L^2_\loc(\RR^d)$, we set
\begin{equation}\label{L2ulnorm}
   \|u\|_{L^2_\ul} \,=\, \sup_{\xi \in \RR^d} \Bigl(\int_{B(\xi,1)} 
   |u(x)|^2  dx \Bigr)^{1/2} \,=\, \sup_{\xi \in \RR^d} 
   \|u\|_{L^2(B(\xi,1))}~.
\end{equation}
The uniformly local Lebesgue space is defined as
\begin{equation}\label{L2uldef}
   L^2_\ul(\RR^d) \,=\, \Bigl\{u \in L^2_\loc(\RR^d) \,\Big|\, 
   \|u\|_{L^2_\ul} < \infty\,,~\lim_{\xi \to 0}
   \|u(\cdot-\xi) - u\|_{L^2_\ul} = 0\Bigr\}~,
\end{equation}
In a similar way, for any $k \in \NN$, we introduce 
the uniformly local Sobolev space
\begin{equation}\label{Hkuldef}
  H^k_\ul(\RR^d) \,=\, \Bigl\{u \in H^k_\loc(\RR^d) \,\Big|\,
  \partial^j_{x_i} u \in L^2_\ul(\RR^d) \hbox{ for } i=1,\ldots,d\text{ and }j
= 0,1,\dots,k \Bigr\}~,
\end{equation}
which is equipped with the natural norm $\|u\|_{H^k_\ul} \,=\,
\Bigl(\sum_{i=1}^d\sum_{j=0}^k \|\partial^j_{x_i} u\|_{L^2_\ul}^2 \Bigr)^{1/2}$.
As shown in \cite{TG-RJ}, the damped wave equation \eqref{eq} is well defined
on $H^1_\ul(\RR^d)\times L^2_\ul(\RR^d)$. The assumption $\lim_{\xi \to 0}
\|u(\cdot-\xi) - u\|_{L^2_\ul} = 0$ in \eqref{L2uldef} introduces a
continuity with respect to translations, which plays the role of the uniform
continuity for continuous functions. It could be possible to work without
this assumption, however $H^1_\ul(\RR^d)$ would not be dense in
$L^2_\ul(\RR^d)$ in this case, which is troublesome. That is why the assumption
$\lim_{\xi \to 0} \|u(\cdot-\xi) - u\|_{L^2_\ul} = 0$ in \eqref{L2uldef} may be
important for the functional analysis. 

We have the following result.
\begin{theorem}
Assume that the assumptions of Theorem \ref{th-main} hold. Then the semigroup
generated by the damped wave equation \eqref{eq} on $H^1_\ul(\RR^d)\times
L^2_\ul(\RR^d)$ is exponentially decreasing. There exist $M$ and
$\lambda>0$ such that, for any solution $U(t)=(u(t),\partial_t u(t))$ of
\eqref{eq} with $U(0)\in H^1_\ul(\RR^d)\times L^2_\ul(\RR^d)$, we have 
$$
\|U(t)\|_{H^1_\ul\times L^2_\ul}\leq Me^{-\lambda t} \|U(0)\|_{H^1_\ul\times
L^2_\ul}
$$
\end{theorem}
\begin{demo}
It is sufficient to show that there exists a time $T>0$ and $C\in
(0,1)$ such that 
\begin{equation}\label{eq-dem-111}
\|U(T)_{|B(\xi,1)}\|_{H^1\times L^2}\leq C \|U(0)\|_{H^1_\ul\times
L^2_\ul}
\end{equation}
for all solutions of \eqref{eq} and all $\xi\in\RR^d$. Due to the finite speed
of propagation of informations and due to the bounds on the metric $K(x)$,
$U(T)_{|B(\xi,1)}$ only depends on the values of $U(0)_{|B(\xi,1+\kappa T)}$
for some $\kappa>0$. Applying Theorem \ref{th-main} to the solution
corresponding to a compactly supported truncation of $U(0)$, we have that 
$$ \|U(T)_{|B(\xi,1)}\|_{H^1\times L^2}\leq Me^{-\lambda T}
\|U(0)_{|B(\xi,2+\kappa T)}\|_{H^1\times L^2}\leq N e^{-\lambda T} T^d
\|U(0)\|_{H^1_\ul\times L^2_\ul}$$
since the ball $B(\xi,2+\kappa T)$ can be covered by a number $\Oc(T^d)$ of
balls of radius $1$. For $T$ large enough, we obtain \eqref{eq-dem-111}, which
shows the result.\end{demo}

\subsection{Linear control}
By HUM Method of Lions (see \cite{Lions}), the exponential decay of the linear
semigroup $e^{At}$ is equivalent to the controllability of the linear wave
equation (in large time). We denote by $\Un_\omega$ the function $\Un_\omega\equiv 1$ on
$\omega$ and $0$ elsewhere. 
\begin{corollary}
Let $\omega$ be an open subset of $\RR^d$ and assume that the hypotheses of
Theorem \ref{th-main} hold with $\gamma=\Un_\omega$. Then, there exists $T>0$
such that, for any $(u_0,u_1)\in
H^1(\RR^d)\times L^2(\RR^d)$ and any $(\tilde u_0,\tilde u_1)\in
H^1(\RR^d)\times L^2(\RR^d)$, there exists a control
$v\in L^1((0,T),L^2(\omega))$ such that the solution $u$ of
$$
\left\{\begin{array}{ll}
\partial^2_{tt} u - \div(K(x)\grad u)+u=\Un_\omega v(x,t)
~~~~&(t,x)\in (0,T) \times\RR^d~,\\
(u,\partial_t u)(\cdot,0)=(u_0,u_1)&
\end{array}\right.
$$
satisfies $(u,\partial u)(\cdot,T)=(\tilde u_0,\tilde u_1)$.
\end{corollary}

\subsection{Global attractor and stabilisation for the non-linear
equation}
Another related problem is the asymptotic behaviour of the non-linear equation
as studied in \cite{Zuazua2}, \cite{Feireisl}, \cite{DLZ} or \cite{Jol-Lau}. One
considers the
non-linear equation
\begin{equation}\label{eq-nl}
\left\{\begin{array}{ll}
\partial^2_{tt} u +\gamma(x) \partial_t u=
\div(K(x)\grad u)-u - f(x,u) ~~~~&(t,x)\in\RR_+\times\RR^d~,\\
(u,\partial_t u)(\cdot,0)=(u_0,u_1)\in H^1(\RR^d)\times L^2(\RR^d)&
\end{array}\right.
\end{equation}
with $f(x,s)\in\Cc^1(\RR^d\times\RR,\RR)$ compactly supported in
$x$, satisfying
\begin{equation}\label{hyp-f-1}
|f(x,s)|\leq C(1+|s|)^p~\text{ and }~~|f'(x,s)|\leq C(1+|s|)^{p-1}
\end{equation}
with $1\leq p <(d+2)/(d-2)$ (or any $p\geq 1$ if $d<3$) and 
\begin{equation}\label{hyp-f-2}
\liminf_{|s|\rightarrow +\infty} ~\max_{x\in\supp(f)}~ f(x,s)s~\geq~ 0~.
\end{equation}
To each solution of \eqref{eq-nl}, one can associate the energy
$$E(u):=E(u,\partial_t u)=\frac 12  \int_{\RR^d}(|\partial_t u|^2+|\grad
u\t.K(x).\grad u|+|u|^2) + \int_{\RR^d}V(x,u)~,$$
where $V(x,u)=\int_0^u f(x,s)ds$. This energy is non-increasing since
$$\partial_t E(u(t))=-\int_{\RR^d} \gamma(x)|\partial_t u(x,t)|^2\,dx~.$$
Notice that \eqref{hyp-f-1} implies that the energy is well defined and
that \eqref{hyp-f-2} shows that the energy is bounded from below and that the
bounded sets of $X=H^1(\RR^d)\times L^2(\RR^d)$ are equivalent to the sets of
bounded energy. Thus any trajectory of \eqref{eq-nl} has a non-increasing
energy and stays bounded in $X$. Assume now that the semigroup $e^{At}$ is
exponentially stable, then any trajectory $U=(u,\partial_t u)$ satisfies
\begin{equation}\label{eq-inte}
U(t)=e^{At}U(0)+\int_0^t
e^{A(t-s)}\left(\begin{array}{c}0\\f(x,u(x,s))\end{array}\right)\,ds~. 
\end{equation}
The first term of \eqref{eq-inte} is decaying exponentially fast and the
integral term is compact since $f$ is compactly supported in $x$ and due to
either the compact Sobolev embedding $H^1\hookrightarrow L^{2p}$ for
$p<d/(d-2)$ or to more technical arguments based on the Strichartz estimates
for $p\in[d/(d-2),(d+2)/(d-2))$ (see \cite{DLZ} and see \cite{Jol-Lau}). Thus,
following the ideas of
\cite{DLZ} and \cite{Jol-Lau}, we obtain that any solution is asymptotically
compact and converges to a trajectory with constant energy. Now, we would like
to show that the energy $E$ associated to \eqref{eq-nl} is a Lyapounov function,
that is that it is non-increasing and cannot be constant along a solution
$u(t)$, except of course if $u(t)$ is an equilibrium point. If \eqref{eq-nl}
admits a Lyapounov function, one says that the corresponding dynamical system is
gradient. In particular, it cannot admit periodic orbits, homoclinic
orbits\ldots{} The gradient structure of \eqref{eq-nl}, together with its
asymptotic compactness, will also ensure the
existence of a compact global attractor, that is a compact invariant set of $X$
which attracts all the trajectories of \eqref{eq-nl}. This set is a central
object of the theory of dynamical systems. It contains all the solutions $u(t)$,
which exist for all $t\in\RR$ and which are uniformly bounded in
$H^1(\RR^d)\times L^2(\RR^d)$ for all $t\in\RR$ (as equilibrium points,
heteroclinic orbits etc.). See for example \cite{Hale-book} and \cite{Raugel}
for a review on the concepts of compact global attractors, of asymptotic
compactness or of gradient structure.

To show that the energy $E$ is a Lyapounov function, that is that it cannot be
constant along a trajectory, except if this trajectory is an equilibrium point,
one has to use a suitable unique continuation property.
In \cite{Zuazua2} and \cite{DLZ}, the authors use a unique continuation
property, which needs geometric assumptions stronger than the one required for
the exponential decay of $e^{At}$. However, we have shown in \cite{Jol-Lau} that
the geometric assumptions required for the exponential decay of $e^{At}$ are
sufficient if we assume that $f$ is smooth and partially analytic. Thus, we can
improve the result of \cite{Jol-Lau} by using a weaker assumption than the one
that $\gamma\geq \alpha>0$ outside a compact set.
\begin{corollary}\label{coro-1}
Assume that the hypotheses of Theorem \ref{th-main} hold. Also assume that
$f\in\Cc^\infty(\RR\times\RR,\RR)$ satisfies \eqref{hyp-f-1} and \eqref{hyp-f-2}
and that $f$ is compactly supported in $x$ and analytic with respect to $u$.
Then the dynamical system generated by \eqref{eq-nl} in $H^1_0(\RR^d)\times
L^2(\RR^d)$ is gradient and admits a compact global attractor $\Ac$.

Moreover, if $f(x,u)u\geq 0$ for any $(x,u)\in\RR^{d+1}$, then the semilinear
damped wave equation \eqref{eq-nl} is stabilised in the sense that for any
$E_0\geq 0$, there exist $K>0$ and $\lambda>0$ such
that, for all solutions $u$ of \eqref{eq-nl} with $E(u(0))\leq E_0$, 
$E(u(t))\leq Me^{-\lambda t} E(u(0))$ for all $t\geq 0$.
\end{corollary}

Of course, in the one-dimensional case $d=1$, the unique continuation property
holds without any additional assumption. We obtain the following result, where
one can omit the assumption $\gamma\geq \alpha>0$ close to $\pm\infty$ used in
\cite{Zuazua2}.
\begin{corollary}\label{coro-2}
Assume that the hypotheses of Theorem \ref{th-main} hold in dimension $d=1$ and
that $f\in\Cc^1(\RR)$ satisfies \eqref{hyp-f-1} and \eqref{hyp-f-2}
and that $f$ is compactly supported in $x$. Then the dynamical system generated
by \eqref{eq-nl} in $H^1_0(\RR)\times L^2(\RR)$ is gradient and admits a
compact global attractor $\Ac$.

Moreover, if $f(x,u)u\geq 0$ for any $(x,u)\in\RR^{2}$, then the semilinear
damped wave equation \eqref{eq-nl} is stabilised.
\end{corollary}

Notice that Corollaries \ref{coro-1} and \ref{coro-2} are not
exactly generalisations of
the results of \cite{Jol-Lau} and \cite{Zuazua2}. Indeed, $f$ is assumed to be
compactly supported in $x$. In \cite{Jol-Lau} and \cite{Zuazua2}, because of the
assumption $\gamma \geq \alpha>0$ outside a compact set, one is able to deal
with nonlinearities $f$ being not compactly supported.

This type of non-linear stabilisation results is also closely related to the
problem of global control of the non-linear wave equation, see \cite{DLZ}, 
\cite{Jol-Lau} and \cite{Jol-Lau-2}. For example, one gets the following result
in dimension $d=1$.
\begin{corollary}
Let $\omega$ be an open subset of $\RR$. Assume that there exist $L>0$ and
$\varepsilon>0$ such that $\omega$ contains an interval of length $\varepsilon$
in any interval $[x,x+L]$, $x\in\RR$. Let also $f\in\Cc^1(\RR\times\RR,\RR)$
compactly supported in $x$ and satisfying \eqref{hyp-f-2}. 

Then, for all $E_0\geq 0$, there exists $T>0$ such that, for any $(u_0,u_1)$
and $(\tilde u_0,\tilde u_1)\in H^1(\RR)\times L^2(\RR)$ with energy $E$
less than $E_0$, there exists a control $v\in L^1((0,T),L^2(\omega))$ such that
the solution $u$ of
$$
\left\{\begin{array}{ll}
\partial^2_{tt} u - \div(K(x)\grad u)+u+f(x,u)=\Un_\omega v(x,t)
~~~~&(t,x)\in (0,T) \times\RR~,\\
(u,\partial_t u)(\cdot,0)=(u_0,u_1)&
\end{array}\right.
$$
satisfies $(u,\partial u)(\cdot,T)=(\tilde u_0,\tilde u_1)$.
\end{corollary}

%%%%%%%%%%%%%%%%%%%%%%%%%%%%%%%%%%%%%%%%%%%%
%%%%%%%%%%%%%%%%%%%%%%%%%%%%%%%%%%%%%%%%%%%%

\setcounter{section}{0}
\renewcommand{\thesection}{{\Alph{section}}}

\section{Appendix: pseudo-differential semiclassical
calculus}\label{sect-A-pseudo}

In section, we recall the main results and notations of pseudo-differential
calculus, which are used in this paper. The details and the proofs could be
found in many textbooks, as \cite{Zworski}, \cite{Martinez}, \cite{Alinhac}
or \cite{Lerousseau}.

Let $h>0$ be a small parameter, say that $h\in(0,1]$.
We say that $a(x,\xi)\in\Cc^\infty(\RR^d\times\RR^d)$ is a {\it symbol} of order
$m$ if, for any multi-indices $\alpha$ and $\beta$, there exists
$C_{\alpha,\beta}$ such that 
$$\sup_{(x,\xi)\in\RR^{2d}} |\partial_x^\alpha\partial_\xi^\beta
a(x,\xi)|\leq C_ {\alpha,\beta} (1+|\xi|)^{m-|\beta|}~$$
and $m$ is the smallest number such that this bounds holds. To each symbol $a$,
we associate the {\it pseudo-differential semiclassical operator} denoted by
$\Op_h(a)$ and defined by Weyl quantization
\begin{equation}\label{eq-weyl}
\Op_h(a)u=\frac{1}{(2\pi)^n}
\int_{\RR^{2d}}e^{i(x-y).\xi}\,a\left(\frac{x+y}2,h\xi\right)u(y)\,dy\,d\xi~.
\end{equation}
If $a$ is of order $m$, then, for any $s>0$, $\Op_h(a)$ is a bounded
operator from $H^s(\RR^d)$ into $H^{s-m}(\RR^d)$, uniformly with respect to
$h\in(0,1]$. 

Let $f\in\Cc^\infty(\RR^d,\RR)$ be a smooth bounded function such that all its
derivative are also bounded functions of $L^\infty(\RR^d)$. It is not so trivial
but well known that the simple operator $u\mapsto f(x)u$  has
for symbol $f(x)$, which is of order $0$ (see for example Chapter 4 of
\cite{Zworski}). More classically, we have that the operator $h\grad$ has for
symbol $i\xi$, which is of order $1$. Using Proposition \ref{prop-compo} below,
one can check that the operator $h^2\Delta_K=h^2\div(K(x)\grad\cdot)$ has for
principal symbol $-\xi\t.K(x).\xi$ in the sense that
$$h^2\Delta_K=\Op_h(-\xi\t.K(x).\xi)+\mathcal{O}_{L^2\rightarrow L^2}(h^2)~.$$

Composing two pseudo-differential operators, one obtain a pseudo-differential
operator, which symbol can be express by an asymptotic development. In this
paper, we will simply use the following cases, see \cite{Zworski},
\cite{Martinez}, \cite{Alinhac} or any other textbooks on pseudo-differential
calculus for more precise developments and for proofs.
\begin{prop}[Composition]\label{prop-compo}
Let $a$ and $b$ be two symbols of order $m$ and $n$ respectively. Assume
that $m+n\leq 2$, then $\Op_h(a)\circ \Op_h(b)$ is a pseudo-differential
operator of order $m+n$ and its symbol $(a\# b)$ satisfies
$$ a\#b= ab-\frac{ih}2\{a,b\}+\mathcal{O}_{L^2\rightarrow L^2}(h^2)$$
where $\{a,b\}=\partial_\xi a\partial_x b-\partial_\xi b\partial_x a$ is the
Poisson bracket of $a$ and $b$ and is of order $m+n-1$. In particular, if
$m+n\leq 1$, then
$$\Op_h(a)\circ \Op_h(b)=\Op_h(ab)+\mathcal{O}_{L^2\rightarrow L^2}(h)~.$$
\end{prop}

\begin{corollary}[Commutators]\label{coro-commut}$~$
\begin{itemize}
\item[i)] If $a$ is of order $1$ or less and if $b$ is of order $0$, then the
commutator $[\Op_h(a),\Op_h(b)]=\Op_h(a)\circ
\Op_h(b)-\Op_h(b)\circ \Op_h(a)$ is of order $0$ or less and of estimate
$\mathcal{O}_{L^2\rightarrow L^2}(h)$.
\item[ii)] If $a$ is of order $2$ and if $b$ is of order $0$, then their
commutator
is of order $1$ and 
$$[\Op_h(a),\Op_h(b)]=-ih\Op_h(\{a,b\})+ \mathcal{O}_{L^2\rightarrow
L^2}(h^2)~.$$
\end{itemize} 
\end{corollary}

\begin{corollary}[Inverse]\label{coro-inverse}                                            
Assume that $a$ is a symbol of order $m\geq 0$ such that there exists
$\alpha>0$ such that $|a(x,\xi)|\geq\alpha$ for all $(x,\xi)\in\RR^{2d}$. Then
$b=1/a$ is a symbol of order $-m$ and $\Op_h(b)$ is a first order
inverse of $a$ in the sense that 
$$\Op_h(a)\circ\Op_h(b)=Id+\mathcal{O}_{L^2\rightarrow L^2}(h)~~~\text{
and }~~~
\Op_h(b)\circ\Op_h(a)=Id+\mathcal{O}_{L^2\rightarrow L^2}(h)
$$
In particular, $\Op_h(a)$ is invertible for $h$ sufficiently small.
\end{corollary}

Simply using the definition \eqref{eq-weyl} and a straightforward computation,
we obtain the expression of the adjoint operator.
\begin{prop}[Adjoint operator]\label{prop-adjoint}
Let $a$ be a symbol of order $m\geq 0$ and $u\in H^m(\RR^d)$. Then,
$$\left(\Op_h(a)\right)^*=\Op_h(\overline{a})~~~~\text{ and }~~~~\Re(\langle
\,\Op_h(a) u\,|\,u\,\rangle_{L^2})=\langle\,
\Op_h(\Re(a))u\,|\,u\,\rangle_{L^2}~.$$
\end{prop}

We will also use a version of G\aa{}rding inequality. We give a
complete proof since in many cases, the result is stated for functions defined
in a compact domain rather than in $\RR^d$, which avoids the question of
uniformity of the constants.
\begin{prop}[G\aa{}rding inequality]\label{prop-garding}
Let $a$ be a symbol of order $m\geq 0$ such that there exists a positive
constant $\alpha$ and $0\leq k\leq m$ such that, for all $(x,\xi)\in\RR^{2d}$,
$\Re(a(x,\xi))\geq \alpha(1+|\xi|)^{k} >0$. Then, there exists $c>0$ such
that, for any $u\in H^m(\RR^d)$ and any $h>0$ sufficiently small, 
$$\Re(\langle \,\Op_h(a) u\,|\,u\,\rangle_{L^2})~\geq ~
c\,\left(\|u\|_{L^2}^2+h^{k/2}\|u\|_{H^{k/2}}^2\right)~.$$
\end{prop}
\begin{demo}
Assume first that $k=0$. We define $b(x,\xi)= \sqrt{\Re(a(x,\xi))}$. We notice
that $b$ is a well defined symbol of order $m/2$ and $\Op_h(b)$ is invertible in
the sense of Corollary \ref{coro-inverse}. In particular, there exists
$\kappa>0$
such that $\|\Op_h(b)u\|\geq \kappa \|u\|^2$ for any $h$ small enough. Using
Proposition \ref{prop-adjoint}, we get that
\begin{align*}
\Re(\langle \,\Op_h(a) u\,|\,u\,\rangle_{L^2})& = \langle \,\Op_h(b^2)
u\,|\,u\,\rangle_{L^2} = \langle \,\Op_h(b)^2 u\,|\,u\,\rangle_{L^2} +
\mathcal{O}(h)\|u\|^2_{L^2}\\ & = \langle \,\Op_h(b)
u\,|\,\Op_h(b)u\,\rangle_{L^2} +
\mathcal{O}(h)\|u\|^2_{L^2}\\
&\geq (\kappa^2+\mathcal{O}(h)) \|u\|^2_{L^2}~,
\end{align*}
which concludes for $h$ small enough.

Now, if $k>0$, we consider $\tilde a=a-\beta(1+|\xi|^2)^{k/2}$, which
satisfies the proposition for $k=0$ if $\beta$ is small enough. We get that 
$$ \Re(\langle \,\Op_h(a) u\,|\,u\,\rangle_{L^2})~-~ \beta \langle
\,\Op_h((1+|\xi|^2)^{k/2}) u\,|\,u\,\rangle_{L^2}~\geq ~
c\|u\|_{L^2}^2~.$$
To conclude, we only have to remark that $\langle
\,\Op_h((1+|\xi|^2)^{k/2}) u\,|\,u\,\rangle_{L^2}$ is equivalent to
$\|u\|_{L^2}^2+h^{k/2}\|u\|_{H^{k/2}}^2$.
\end{demo}

It is also common to use G\aa{}rding inequality for functions vanishing on some
part of $\RR^d$. In this case, we have to be more careful
about the uniformity of the positive constants, which leads use to work with
$k=m$ .
\begin{prop}[G\aa{}rding inequality with truncation]\label{prop-garding-2}
Let $a$ be a symbol of order $m\geq 0$ and $\omega$ be a subset of $\RR^d$.
Assume that there exists a positive constant $\alpha$ such that, for all
$x\in \RR^{d}\setminus \omega$ and all $\xi\in\RR^d$,
$\Re(a(x,\xi))\geq \alpha(1+|\xi|)^{m}>0$. Then, there exists $c>0$ such
that, for any $u\in H^m(\RR^d)$ such that $u_{|\omega}\equiv 0$ and any $h>0$
sufficiently small, 
$$\Re(\langle \,\Op_h(a) u\,|\,u\,\rangle_{L^2})~\geq ~
c\,\left(\|u\|_{L^2}^2+h^{m/2}\|u\|_{H^{m/2}}^2\right)~.$$
\end{prop}
\begin{demo}
We set
$\Omega_\varepsilon=\{x\in\RR^d,~ d(x,\RR^d\setminus\omega) < \varepsilon \}$.
Notice that, by definition of a symbol of order $m\geq 0$, we have $\partial_x
a(x,\xi)\leq C (1+|\xi|^2)^{m/2}$, which implies that we still have
$\Re(a(x,\xi))\geq \tilde\alpha(1+|\xi|^2)^{m/2}>0$ for $\tilde\alpha
\in(0,\alpha)$ in $\Omega_\varepsilon$ for some small $\varepsilon>0$. Due to
the distance $\varepsilon>0$ between $\RR^d\setminus \Omega_\varepsilon$ and
$\RR^d\setminus \omega$, we can construct a function 
$\chi\in\Cc^\infty_b(\RR^d,[0,1])$ with support in $\omega$ and which is equal
to $1$ outside $\Omega_\varepsilon$. It is then sufficient to apply Proposition
\ref{prop-garding} to $\tilde a(x,\xi)= a(x,\xi)+(1+|\xi|^2)^{m/2}\chi(x)$.
\end{demo}

%%%%%%%%%%%%%%%%%%%%%%%%%%%%%%%%%%%%%%%%%%%%%%%%%%%%%%%%%%%%%%%%%%%%%%%%%%%%
%%%%%%%%%%%%%%%%%%%%%%%%%%%%%%%%%%%%%%%%%%%%%%%%%%%%%%%%%%%%%%%%%%%%%%%%%%%%

\end{document}